\documentclass[11pt]{article}
\usepackage{makeidx}
\usepackage{amsfonts}
\usepackage{amsmath}
\usepackage{amssymb}
\usepackage{graphicx}
\usepackage{bbold}
\usepackage[colorlinks]{hyperref}
\usepackage{palatino}

\hypersetup{colorlinks=true, linkcolor=cyan, citecolor=green,
filecolor=black, urlcolor=black }
\providecommand{\U}[1]{\protect\rule{.1in}{.1in}}
\newtheorem{theorem}{Theorem}

\newtheorem{corollary}[theorem]{Corollary}

\newtheorem{definition}[theorem]{Definition}

\newtheorem{proposition}[theorem]{Proposition}
\newtheorem{remark}[theorem]{Remark}

\newenvironment{proof}[1][Proof]{\noindent\textbf{#1.} }{\ \rule{0.5em}{0.5em}}

\input{tcilatex}

\begin{document}

\title{{Backward stochastic variational inequalities with locally bounded
generators} \thanks{{\scriptsize The work for this paper was supported by
IDEAS project, no. 241/05.10.2011 \smallskip}}}
\author{Lucian Maticiuc$^{{\scriptsize {a,b,\ast}}}$, Aurel R\u{a}\c{s}canu$%
^{{\scriptsize {b,c,\ast}}}$, Adrian Z\u{a}linescu$^{{\scriptsize b,c,\ast}}$
\bigskip \\
$^{{\scriptsize a \,}}${\scriptsize Department of Mathematics,
\textquotedblleft Gheorghe Asachi\textquotedblright\ Technical University of
Ia\c{s}i, Carol I Blvd., no. 11, 700506, Romania,}\\
$^{{\scriptsize b \,}}${\scriptsize Faculty of Mathematics,
\textquotedblleft Alexandru Ioan Cuza\textquotedblright\ University of Ia%
\c{s}i, Carol I Blvd., no. 11, 700506, Romania,}\\
$^{{\scriptsize c \,}}${\scriptsize \textquotedblleft Octav
Mayer\textquotedblright Mathematics Institute of the Romanian Academy, Ia%
\c{s}i, Carol I Blvd., no. 8, 700506, Romania.}}
\maketitle

\begin{abstract}
The paper deals with the existence and uniqueness of the solution of the
backward stochastic variational inequality:%
\begin{equation*}
\left\{
\begin{array}{l}
-dY_{t}+\partial \varphi \left( Y_{t}\right) dt\ni F\left(
t,Y_{t},Z_{t}\right) dt-Z_{t}dB_{t},\;0\leq t<T\medskip  \\
Y_{T}=\eta ,%
\end{array}%
\,\right.
\end{equation*}%
where $F$ satisfies a local boundedness condition.
\end{abstract}

\footnotetext{\textit{\scriptsize E-mail addresses:} {\scriptsize %
lucian.maticiuc@ymail.com (Lucian Maticiuc),~aurel.rascanu@uaic.ro (Aurel R%
\u{a}\c{s}canu),}}

\footnotetext{{\scriptsize adrian.zalinescu@gmail.com (Adrian Z\u{a}linescu).%
}}

\textbf{AMS Classification subjects: }60H10, 93E03, 47J20, 49J40.\medskip


\section{Introduction}

We consider the following backward stochastic variational inequality (BSVI)%
\begin{equation}
\left\{
\begin{array}{l}
-dY_{t}+\partial\varphi\left( Y_{t}\right) dt\ni F\left(
t,Y_{t},Z_{t}\right) dt-Z_{t}dB_{t},\;0\leq t<T\medskip \\
Y_{T}=\eta.%
\end{array}
\right.  \label{BSVI}
\end{equation}
where $\left\{ B_{t}:t\geq0\right\} $ is a standard Brownian motion, $%
\partial\varphi$ is the subdifferential of a convex l.s.c. function $\varphi
$, and $T>0$ is a fixed deterministic time.

The study of the backward stochastic differential equations (equation of
type (\ref{BSVI}) without the subdifferential operator) was initiated by E.
Pardoux and S. Peng in \cite{pa-pe/90} (see also \cite{pa-pe/92}) where is
proved the existence and the uniqueness of the solution for the BSDE under
the assumption of Lipschitz continuity of $F$ with respect to $y$ and $z$
and square integrability of $\eta$ and $F\left( t,0,0\right) $.

The more general case of scalar BSDE with one-sided reflection and
associated optimal control problems was considered by N. El Karoui, C.
Kapoudjian, E. Pardoux, S. Peng, M.C. Quenez in \cite{ka-ka-pa/97} and with
two-sided reflection associated with stochastic game problem by Cvitanic and
Karatzas \cite{cv-ka/96}.

Multidimensional BSDE reflected at the boundary of a convex set was studied
in A. Gegout-Petit and E. Pardoux, \cite{ge-pa/96}.

The standard work on BSVI is that of E. Pardoux and A. R\u{a}\c{s}canu \cite%
{pa-ra/98}, which give a proof of existence and uniqueness of the solution
for (\ref{BSVI}) under the following assumptions on $F:$ monotonicity with
respect to $y$ (in the sense that $\langle y^{\prime }-y,F(t,y^{\prime
},z)-F(t,y,z)\rangle \leq \alpha |y^{\prime }-y|^{2}$), Lipschitzianity with
respect to $z$ and a sublinear growth for $F\left( t,y,0\right) :$%
\begin{equation*}
\left\vert F\left( t,y,0\right) \right\vert \leq \beta _{t}+L\left\vert
y\right\vert ,\;\forall \left( t,y\right) \in \left[ 0,T\right] \times
\mathbb{R}^{m}.
\end{equation*}%
It is proved that there exists a unique triple $\left( Y,Z,K\right) $ such
that%
\begin{equation*}
Y_{t}+K_{T}-K_{t}=\eta +\int_{t}^{T}F\left( s,Y_{s},Z_{s}\right)
ds-\int_{t}^{T}Z_{s}dB_{s},\text{\ a.s., with }dK_{t}\in \partial \varphi
\left( Y_{t}\right) dt.
\end{equation*}%
Moreover the process $K$ is absolute continuous with respect to $dt$. In
\cite{pa-ra/99} the same authors extend the results from \cite{pa-ra/98} to
a Hilbert spaces framework. Using a mixed Euler-Yosida scheme, Maticiuc and
Rotenstein provided in \cite{ma-ro/12} numerical results concerning the multi%
{\footnotesize \-}valued stochastic differential equation (\ref{BSVI}).

Our paper generalize the previous existence and uniqueness results for (\ref%
{BSVI}) by assu\-ming a local boundedness condition (instead of sublinear
growth of $F$), i.e.%
\begin{equation*}
\mathbb{E}\Big(\int _{0}^{T}F_{\rho}^{\#}(s)ds\Big)^{p}<\infty,\text{ where }%
F_{\rho}^{\#}\left( t\right) \overset{def}{=}\sup\limits_{\left\vert
y\right\vert \leq\rho }\left\vert F(t,y,0)\right\vert .
\end{equation*}
Concerning to this requirement on $F$ we remark that a similar one was
considered by E. Pardoux in \cite{pa/99} for the study of BSDE. More
precisely, his result is the following: If $\eta\in L^{2}\left( \Omega;%
\mathbb{R}^{m}\right) $, $F\left( t,0,0\right) \in L^{2}\left( \Omega\times%
\left[ 0,T\right] ;\mathbb{R}^{m}\right) $, $F$ is monotone with respect to $%
y$, Lipschitz with respect to $z$ and there exists a deter\-mi\-nis\-tic
continuous increasing function $\psi$ such that $\forall\left( t,y\right) \in%
\left[ 0,T\right] \times\mathbb{R}^{m},\;\left\vert F\left( t,y,0\right)
\right\vert \leq\left\vert F\left( t,0,0\right) \right\vert +\psi\left(
\left\vert y\right\vert \right) ,\mathbb{P}$-a.s, then there exist a unique
solution for BSDE (\ref{BSVI}) with $\varphi\equiv0$. This result was
generalized by Ph. Briand, B. Delyon, Y. Hu, E. Pardoux, L. Stoica in \cite%
{br-de-hu-pa/03}.

The article is organized as follows: in the next Section we prove some a
priori estimates and the uniqueness result for the solution of BSVI (\ref%
{BSVI}). Section 3 is concerned on the existence result under two
alternative assumptions (which allow to obtain the absolute continuity of
the process $K$) and Section 4 establishes the general existence result. In
the Appendix we presents, following \cite{pa-ra/08}, some results useful
throughout the paper.

\section{Preliminaries; a priori estimates and the uniqueness result}

Let $\left\{ B_{t}:t\geq0\right\} $ be a $k$-dimensional standard Brownian
motion defined on some complete probability space $(\Omega,\mathcal{F},%
\mathbb{P})$. We denote by $\left\{ \mathcal{F}_{t}:t\geq0\right\} $ the
na\-tu\-ral filtration generated by $\left\{ B_{t}:t\geq0\right\} $\ and
augmented by $\mathcal{N}$, the set of$\;\mathbb{P}$- null events of $%
\mathcal{F}$:%
\begin{equation*}
\mathcal{F}_{t}=\sigma\{B_{r}:0\leq r\leq t\}\vee\mathcal{N}.
\end{equation*}
We suppose that the following assumption holds

\begin{itemize}
\item[\textrm{(A}$_{1}$\textrm{)}] $\eta:\Omega\rightarrow\mathbb{R}^{m}$ is
a $\mathcal{F}_{T}$-measurable random vector,

\item[\textrm{(A}$_{2}$\textrm{)}] $F:\Omega \times \left[ 0,T\right] \times
\mathbb{R}^{m}\times \mathbb{R}^{m\times k}\rightarrow \mathbb{R}^{m}$
satisfies that, for all $y\in \mathbb{R}^{m}$, $z\in \mathbb{R}^{m\times k}$%
, $\left( \omega ,t\right) \longmapsto F\left( \cdot ,\cdot ,y,z\right)
:\Omega \times \left[ 0,T\right] \rightarrow \mathbb{R}^{m}$ is
progressively measurable stochastic process, and there exist $\mu :\Omega
\times \left[ 0,T\right] \rightarrow \mathbb{R}$ and $\ell :\Omega \times %
\left[ 0,T\right] \rightarrow \mathbb{R}_{+}$ progressively measurable
stochastic processes with%
\begin{equation*}
\int_{0}^{T}\left( \left\vert \mu _{t}\right\vert +\ell _{t}^{2}\right)
dt<\infty ,
\end{equation*}%
such that, for all $t\in \left[ 0,T\right] $, $y,y^{\prime }\in \mathbb{R}%
^{m}$ and$\;z,z^{\prime }\in \mathbb{R}^{m\times k},\;\mathbb{P}$-a.s.:%
\begin{equation*}
\begin{array}{cl}
\left( C_{y}\right) & y\longmapsto F\left( t,y,z\right) :\mathbb{R}%
^{m}\rightarrow \mathbb{R}^{m}\text{ is continuous,}\medskip \\
\left( M_{y}\right) & \left\langle y^{\prime }-y,F(t,y^{\prime
},z)-F(t,y,z)\right\rangle \leq \mu _{t}|y^{\prime }-y|^{2},\medskip \\
\left( L_{z}\right) & |F(t,y,z^{\prime })-F(t,y,z)|\leq \ell _{t}|z^{\prime
}-z|,\medskip \\
\left( B_{y}\right) & \displaystyle\int_{0}^{T}F_{\rho }^{\#}\left( s\right)
ds<\infty ,\;\forall \;\rho \geq 0,\medskip \\
& \;\;\;\;\;\;\;\text{where, for }\rho \geq 0,\;F_{\rho }^{\#}\left(
t\right) \overset{def}{=}\sup\limits_{\left\vert y\right\vert \leq \rho
}\left\vert F(t,y,0)\right\vert ,%
\end{array}%
\end{equation*}

\item[\textrm{(A}$_{3}$\textrm{)}] $\varphi :\mathbb{R}^{m}\rightarrow
(-\infty ,+\infty ]$ is a proper, convex l.s.c. function.
\end{itemize}

The subdifferential of $\varphi$ is given by
\begin{equation*}
\partial\varphi\left( y\right) =\left\{ \hat{y}\in\mathbb{R}%
^{m}:\left\langle \hat{y},v-y\right\rangle +\varphi\left( y\right) \leq
\varphi\left( v\right) ,\;\forall~v\in\mathbb{R}^{m}\right\}.
\end{equation*}

We define%
\begin{equation*}
\begin{array}{l}
\mathrm{Dom}\left( \varphi \right) =\left\{ y\in \mathbb{R}^{m}:\varphi
\left( y\right) <\infty \right\} ,\medskip \\
\mathrm{Dom}\left( \partial \varphi \right) =\left\{ y\in \mathbb{R}%
^{m}:\partial \varphi \left( y\right) \neq \emptyset \right\} \subset
\mathrm{Dom}\left( \varphi \right)%
\end{array}%
\end{equation*}%
and by $\left( y,\hat{y}\right) \in \partial \varphi $ we understand that $%
y\in \mathrm{Dom}\left( \partial \varphi \right) $ and $\hat{y}\in \partial
\varphi \left( y\right) $.

Recall that%
\begin{equation*}
\overline{\mathrm{Dom}\left( \varphi \right) }=\overline{\mathrm{Dom}\left(
\partial \varphi \right) },\;\;\mathrm{Int}\left( \mathrm{Dom}\left( \varphi
\right) \right) =\mathrm{Int}\left( \mathrm{Dom}\left( \partial \varphi
\right) \right) .
\end{equation*}%
Let $\varepsilon >0$ and the Yosida regularization of $\varphi :$%
\begin{equation}
\varphi _{\varepsilon }\left( y\right) \overset{def}{=}\inf \left\{ \frac{1}{%
2\varepsilon }\left\vert y-v\right\vert ^{2}+\varphi \left( v\right) :v\in
\mathbb{R}^{m}\right\} =\frac{1}{2\varepsilon }\left\vert y-J_{\varepsilon
}\left( y\right) \right\vert ^{2}+\varphi \left( J_{\varepsilon }\left(
y\right) \right) ,  \label{def Yosida regul}
\end{equation}%
where $J_{\varepsilon }\left( y\right) =\left( I_{m\times m}+\varepsilon
\partial \varphi \right) ^{-1}\left( y\right) $. Remark that $\varphi
_{\varepsilon }$ is a $C^{1}$ convex function and $J_{\varepsilon }$ is a $1$%
-Lipschitz function.

We mention some properties (see H. Br\'{e}zis \cite{br/73}, and E. Pardoux,
A. R\u{a}\c{s}canu \cite{pa-ra/98} for the last one): for all $x,y\in
\mathbb{R}^{m}$%
\begin{equation}
\begin{array}{cl}
(a) & \nabla \varphi _{\varepsilon }(y)=\partial \varphi _{\varepsilon
}\left( y\right) =\dfrac{y-J_{\varepsilon }\left( y\right) }{\varepsilon }%
\in \partial \varphi (J_{\varepsilon }y),\medskip \\
(b) & \left\vert \nabla \varphi _{\varepsilon }(x)-\nabla \varphi
_{\varepsilon }(y)\right\vert \leq \dfrac{1}{\varepsilon }\left\vert
x-y\right\vert ,\medskip \\
(c) & \left\langle \nabla \varphi _{\varepsilon }(x)-\nabla \varphi
_{\varepsilon }(y),x-y\right\rangle \geq 0,\medskip \\
(d) & \left\langle \nabla \varphi _{\varepsilon }(x)-\nabla \varphi _{\delta
}(y),x-y\right\rangle \geq -(\varepsilon +\delta )\left\langle \nabla
\varphi _{\varepsilon }(x),\nabla \varphi _{\delta }(y)\right\rangle%
\end{array}
\label{ineq Yosida}
\end{equation}%
We denote by $\mathcal{S}_{m}^{p}[0,T]$ the space of (equivalent classes of)
progressively measurable and continuous stochastic processes $X:\Omega
\times \left[ 0,T\right] \rightarrow \mathbb{R}^{m}$ such that%
\begin{equation*}
\mathbb{E}\sup_{t\in \left[ 0,T\right] }\left\vert X_{t}\right\vert
^{p}<\infty ,\;\text{if }p>0,
\end{equation*}%
and by $\Lambda _{m}^{p}\left( 0,T\right) $ the space of (equivalent classes
of) progressively measurable stochastic process $X:\Omega \times \left[ 0,T%
\right] \rightarrow \mathbb{R}^{m}$ such that%
\begin{equation*}
\begin{array}{ll}
\displaystyle\int_{0}^{T}\left\vert X_{t}\right\vert ^{2}dt<\infty ,\;%
\mathbb{P}\text{-a.s. }\omega \in \Omega , & \text{if }p=0,\medskip \\
\displaystyle\mathbb{E}\left( \int_{0}^{T}\left\vert X_{t}\right\vert
^{2}dt\right) ^{p/2}<\infty , & \text{if }p>0%
\end{array}%
\end{equation*}%
For a function $g:\left[ 0,T\right] \rightarrow \mathbb{R}^{m}$, let us
denote by $\left\updownarrow g\right\updownarrow _{T}$ the total variation
of $g$ on $\left[ 0,T\right] $ i.e.%
\begin{equation*}
\left\updownarrow g\right\updownarrow _{T}\overset{def}{=}\sup \left\{
\sum\limits_{i=0}^{n-1}\left\vert g\left( t_{i+1}\right) -g\left(
t_{i}\right) \right\vert :n\in \mathbb{N}^{\ast },\;0=t_{0}<t_{1}<\cdots
t_{n}=T\right\} ,
\end{equation*}%
and by $BV\left( \left[ 0,T\right] ;\mathbb{R}^{m}\right) $ the space of the
functions $g:\left[ 0,T\right] \rightarrow \mathbb{R}^{m}$ such that $%
\left\updownarrow g\right\updownarrow _{T}<\infty $ ($BV\left( \left[ 0,T%
\right] ;\mathbb{R}^{m}\right) $ equipped with the norm $\left\vert
\left\vert g\right\vert \right\vert _{BV\left( \left[ 0,T\right] ;\mathbb{R}%
^{m}\right) }\overset{def}{=}\left\vert g(0)\right\vert +\left\updownarrow
g\right\updownarrow _{T}$ is a Banach space).

\begin{definition}
\label{def 1 - BSVI}A pair $\left( Y,Z\right) \in S_{m}^{0}\left[ 0,T\right]
\times \Lambda _{m\times k}^{0}\left( 0,T\right) $ of stochastic processes
is a solution of backward stochastic variational inequality (\ref{BSVI}) if
there exists $K\in S_{m}^{0}\left[ 0,T\right] $ with $K_{0}=0$, such that%
\begin{equation*}
\begin{array}{ll}
\left( a\right) & \displaystyle\left\updownarrow K\right\updownarrow
_{T}+\int_{0}^{T}\left\vert \varphi \left( Y_{t}\right) \right\vert
dt+\int_{0}^{T}\left\vert F(t,Y_{t},Z_{t})\right\vert dt<\infty ,\;\text{a.s.%
},\medskip \\
\left( b\right) & dK_{t}\in \partial \varphi \left( Y_{t}\right) dt,\;\text{%
a.s.}\;\text{that is:}\;\mathbb{P}\text{-a.s.},\medskip \\
\multicolumn{1}{r}{} & \multicolumn{1}{r}{\displaystyle\int\nolimits_{t}^{s}%
\langle y(r)-Y_{r},dK_{r}\rangle +\int\nolimits_{t}^{s}\varphi (Y_{r})dr\leq
\int\nolimits_{t}^{s}\varphi (y(r))dr,\medskip} \\
\multicolumn{1}{r}{} & \multicolumn{1}{r}{\forall y\in C([0,T];\mathbb{R}%
^{d}),\;\forall 0\leq t\leq s\leq T,}%
\end{array}%
\end{equation*}%
and, $\mathbb{P}$-a.s., for all\ $t\in \left[ 0,T\right] :$%
\begin{equation}
Y_{t}+K_{T}-K_{t}=\eta +\int_{t}^{T}F\left( s,Y_{s},Z_{s}\right)
ds-\int_{t}^{T}Z_{s}dB_{s}  \label{BSVI explicit form}
\end{equation}%
(we also say that triplet $\left( Y,Z,K\right) $ is solution of equation (%
\ref{BSVI})).
\end{definition}

\begin{remark}
If $K$ is absolute continuous with respect to $dt$, i.e. there exists a
progre\-ssi\-vely measurable stochastic process $U$ such that%
\begin{equation*}
\int _{0}^{T}\left\vert U_{t}\right\vert dt<\infty,\;\text{a.s. and }%
K_{t}=\int _{0}^{t}U_{s}ds,\;\text{for all }t\in\left[ 0,T\right] ,
\end{equation*}
then $dK_{t}\in\partial\varphi\left( Y_{t}\right) dt$ means
\begin{equation*}
U_{t}\in\partial\varphi\left( Y_{t}\right) ,\text{\ }dt\text{-a.e., a.s.}
\end{equation*}
\end{remark}

If $dK_{t}\in \partial \varphi \left( Y_{t}\right) dt$ and $d\tilde{K}%
_{t}\in \partial \varphi (\tilde{Y}_{t})dt$ then we clearly have%
\begin{equation*}
\int_{0}^{T}\left\vert \varphi \left( Y_{t}\right) \right\vert
dt+\int_{0}^{T}|\varphi (\tilde{Y}_{t})|dt<\infty ,\text{\ a.s.}
\end{equation*}%
and, using the subdifferential inequalities%
\begin{equation*}
\begin{array}{l}
\displaystyle\int\nolimits_{t}^{s}\langle \tilde{Y}_{r}-Y_{r},dK_{r}\rangle
+\int\nolimits_{t}^{s}\varphi (Y_{r})dr\leq \int\nolimits_{t}^{s}\varphi (%
\tilde{Y}_{r})dr,\medskip \\
\displaystyle\int\nolimits_{t}^{s}\langle Y_{r}-\tilde{Y}_{r},d\tilde{K}%
_{r}\rangle +\int\nolimits_{t}^{s}\varphi (\tilde{Y}_{r})dr\leq
\int\nolimits_{t}^{s}\varphi (Y_{r})dr,%
\end{array}%
\end{equation*}%
we infer that, for all $0\leq t\leq s\leq T$%
\begin{equation}
\int\nolimits_{t}^{s}\langle Y_{r}-\tilde{Y}_{r},dK_{r}-d\tilde{K}%
_{r}\rangle \geq 0,\ \text{a.s.}  \label{monoton}
\end{equation}%
Let $a,p>1$ and%
\begin{equation}
V_{t}=V_{t}^{a,p}\overset{def}{=}\int_{0}^{t}\big(\mu _{s}+\dfrac{a}{2n_{p}}%
\ell _{s}^{2}\big)ds,  \label{def V}
\end{equation}%
where $n_{p}=\left( p-1\right) \wedge 1$.

Denote%
\begin{equation*}
S_{m}^{1^{+},p}\left[ 0,T\right] \overset{def}{=}\left\{ Y\in S_{m}^{0}\left[
0,T\right] :\exists~a>1,\ \mathbb{E}\sup\limits_{s\in\left[ 0,T\right]
}|e^{V_{s}^{a,p}}Y_{s}|^{p}<\infty\right\} .
\end{equation*}
Remark that if $\mu_{s}$ and $\ell_{s}^{2}$ are deterministic functions
then, for all $p>1$, $S_{m}^{1^{+},p}\left[ 0,T\right] =S_{m}^{p}\left[ 0,T%
\right] $.

\begin{proposition}
\label{prop 1}Let $\left( u_{0},\hat{u}_{0}\right) \in \partial \varphi $
and assumptions \textrm{(A}$_{1}-$\textrm{A}$_{3}$\textrm{)} be satisfied.
Then for every $a,p>1$ there exists a constant $C_{a,p}$ such that for every
$\left( Y,Z\right) $ solution of BSDE (\ref{BSVI}) satisfying%
\begin{equation*}
\mathbb{E}\sup\limits_{s\in \left[ 0,T\right] }e^{pV_{s}}\left\vert
Y_{s}-u_{0}\right\vert ^{p}+\mathbb{E}\Big(\int_{0}^{T}e^{V_{s}}\left(
\left\vert \hat{u}_{0}\right\vert +\left\vert F\left( s,u_{0},0\right)
\right\vert \right) ds\Big)^{p}<\infty ,
\end{equation*}%
the following inequality holds $\mathbb{P}$-a.s., for all $t\in \left[ 0,T%
\right] :$%
\begin{equation}
\begin{array}{l}
\displaystyle\mathbb{E}^{\mathcal{F}_{t}}\left[ \sup\limits_{s\in \left[ t,T%
\right] }\left\vert e^{V_{s}}\left( Y_{s}-u_{0}\right) \right\vert
^{p}+\left( \int_{t}^{T}e^{2V_{s}}\left\vert Z_{s}\right\vert ^{2}ds\right)
^{p/2}\right] \medskip \\
\displaystyle\;\;+\mathbb{E}^{\mathcal{F}_{t}}\left(
\int_{t}^{T}e^{2V_{s}}\left\vert \varphi (Y_{s})-\varphi \left( u_{0}\right)
\right\vert ds\right) ^{p/2}+\mathbb{E}^{\mathcal{F}_{t}}%
\int_{t}^{T}e^{pV_{s}}\left\vert Y_{s}-u_{0}\right\vert ^{p-2}\mathbf{1}%
_{Y_{s}\neq u_{0}}\left\vert Z_{s}\right\vert ^{2}ds\medskip \\
\displaystyle\;\;+\mathbb{E}^{\mathcal{F}_{t}}\int_{t}^{T}e^{pV_{s}}\left%
\vert Y_{s}-u_{0}\right\vert ^{p-2}\mathbf{1}_{Y_{s}\neq u_{0}}\left\vert
\varphi (Y_{s})-\varphi \left( u_{0}\right) \right\vert ds\medskip \\
\displaystyle\leq C_{a,p}\mathbb{E}^{\mathcal{F}_{t}}\bigg[%
e^{pV_{T}}\left\vert \eta -u_{0}\right\vert ^{p}+\left(
\int_{t}^{T}e^{V_{s}}\left\vert \hat{u}_{0}\right\vert ds\right) ^{p}+\left(
\int_{t}^{T}e^{V_{s}}\left\vert F\left( s,u_{0},0\right) \right\vert
ds\right) ^{p}\bigg]%
\end{array}
\label{ineq prop 1}
\end{equation}%
and, for every $R_{0}>0$ and $p\geq 2,$%
\begin{equation}
\begin{array}{l}
\displaystyle R_{0}^{p/2}\mathbb{E}^{\mathcal{F}_{t}}\left(
\int_{t}^{T}e^{2V_{s}}\left\vert F\left( s,Y_{s},Z_{s}\right) \right\vert
ds\right) ^{p/2}+\mathbb{E}^{\mathcal{F}_{t}}\int_{t}^{T}e^{pV_{s}}\left%
\vert Y_{s}-u_{0}\right\vert ^{p-2}\mathbf{1}_{Y_{s}\neq u_{0}}\left\vert
F\left( s,Y_{s},Z_{s}\right) \right\vert ds\medskip \\
\displaystyle\leq C_{a,p}~\bigg[\mathbb{E}^{\mathcal{F}_{t}}e^{pV_{T}}\left%
\vert \eta -u_{0}\right\vert ^{p}+R_{0}^{p/2}\mathbb{E}^{\mathcal{F}%
_{t}}\left( \int_{t}^{T}e^{2V_{s}}\mathbf{1}_{p\geq 2}\left(
F_{u_{0},R_{0}}^{\#}\left( s\right) +R_{0}\gamma _{s}^{+}\right) ds\right)
^{p/2}\medskip \\
\displaystyle\;\;+\mathbb{E}^{\mathcal{F}_{t}}\left(
\int_{t}^{T}e^{V_{s}}\left( F_{u_{0},R_{0}}^{\#}\left( s\right)
+2R_{0}\left\vert \gamma _{s}\right\vert \right) ds\right) ^{p}\bigg],%
\end{array}
\label{ineq 2 prop 1}
\end{equation}%
where%
\begin{equation*}
F_{u_{0},R_{0}}^{\#}\left( t\right) \overset{def}{=}\sup_{\left\vert
y-u_{0}\right\vert \leq R_{0}}\left\vert F\left( t,y,0\right) \right\vert .
\end{equation*}
\end{proposition}

\begin{proof}
We can write%
\begin{equation*}
Y_{t}-u_{0}=\eta-u_{0}+\int_{t}^{T}\left[ F\left( s,Y_{s},Z_{s}\right)
ds-dK_{s}\right] -\int_{t}^{T}Z_{s}dB_{s}
\end{equation*}
Let $R_{0}\geq0$. The monotonicity property of $F$ implies that, for all $%
\left\vert v\right\vert \leq1:$%
\begin{equation*}
\left\langle F\left( t,u_{0}+R_{0}v,z\right) -F\left( t,y,z\right)
,u_{0}+R_{0}v-y\right\rangle \leq\mu_{t}\left\vert u_{0}+R_{0}v-y\right\vert
^{2},
\end{equation*}
and, consequently%
\begin{equation*}
\begin{array}{l}
R_{0}\left\langle F\left( t,y,z\right) ,-v\right\rangle +\left\langle
F\left( t,y,z\right) ,y-u_{0}\right\rangle \medskip \\
\leq\mu_{t}\left\vert u_{0}+R_{0}v-y\right\vert ^{2}+\left\vert F\left(
t,u_{0}+R_{0}v,z\right) \right\vert \left\vert y-R_{0}v-u_{0}\right\vert
\medskip \\
\leq\mu_{t}\left\vert u_{0}+R_{0}v-y\right\vert ^{2}+\left[
F_{u_{0},R_{0}}^{\#}\left( t\right) +\ell_{t}\left\vert z\right\vert \right]
\left\vert y-R_{0}v-u_{0}\right\vert \medskip \\
\leq\mu_{t}\left\vert u_{0}+R_{0}v-y\right\vert
^{2}+F_{u_{0},R_{0}}^{\#}\left( t\right) \left\vert
y-R_{0}v-u_{0}\right\vert +\dfrac{a}{2n_{p}}\ell_{t}^{2}\left\vert
y-R_{0}v-u_{0}\right\vert ^{2}+\dfrac{n_{p}}{2a}\left\vert z\right\vert
^{2}\medskip \\
\leq F_{u_{0},R_{0}}^{\#}\left( t\right) \left( \left\vert
y-u_{0}\right\vert +R_{0}\right) +\gamma_{t}\left[ \left\vert
y-u_{0}\right\vert ^{2}-2R_{0}\left\langle v,y-u_{0}\right\rangle
+R_{0}^{2}\left\vert v\right\vert ^{2}\right] +\dfrac{n_{p}}{2a}\left\vert
z\right\vert ^{2}\medskip \\
\leq\left[ R_{0}F_{u_{0},R_{0}}^{\#}\left( t\right) +R_{0}^{2}\gamma _{t}^{+}%
\right] +\left[ F_{u_{0},R_{0}}^{\#}\left( t\right) +2R_{0}\left\vert
\gamma_{t}\right\vert \right] \left\vert y-u_{0}\right\vert
+\gamma_{t}\left\vert y-u_{0}\right\vert ^{2}+\dfrac{n_{p}}{2a}\left\vert
z\right\vert ^{2}%
\end{array}%
\end{equation*}
Taking $\sup_{\left\vert v\right\vert \leq1}~,$ we have%
\begin{equation*}
\begin{array}{l}
R_{0}\left\vert F\left( t,Y_{t},Z_{t}\right) \right\vert dt+\left\langle
Y_{t}-u_{0},F\left( t,Y_{t},Z_{r}\right) \right\rangle dt\medskip \\
\leq\left[ R_{0}F_{u_{0},R_{0}}^{\#}\left( t\right) +R_{0}^{2}\gamma _{t}^{+}%
\right] +\left[ F_{u_{0},R_{0}}^{\#}\left( t\right) +2R_{0}\left\vert
\gamma_{t}\right\vert \right] \left\vert Y_{t}-u_{0}\right\vert +\left\vert
Y_{t}-u_{0}\right\vert ^{2}dV_{t}+\dfrac{n_{p}}{2a}\left\vert
Z_{t}\right\vert ^{2}%
\end{array}%
\end{equation*}
From the subdifferential inequalities we have%
\begin{equation*}
\left\vert \varphi(t,Y_{t})-\varphi\left( t,u_{0}\right) \right\vert \leq
\left[ \varphi(t,Y_{t})-\varphi\left( t,u_{0}\right) \right] +2\left\vert
\hat{u}_{0}\right\vert \left\vert Y_{t}-u_{0}\right\vert ,
\end{equation*}
and%
\begin{equation*}
\left[ \varphi(t,Y_{t})-\varphi\left( t,u_{0}\right) \right]
dt\leq\left\langle Y_{t}-u_{0},dK_{t}\right\rangle
\end{equation*}
Therefore%
\begin{equation*}
\left\vert \varphi(t,Y_{t})-\varphi\left( t,u_{0}\right) \right\vert
dt\leq\left\langle Y_{t}-u_{0},dK_{t}\right\rangle +2\left\vert \hat{u}%
_{0}\right\vert \left\vert Y_{t}-u_{0}\right\vert dt.
\end{equation*}
From the above it follows that%
\begin{equation}
\begin{array}{l}
\left[ R_{0}\left\vert F\left( t,Y_{t},Z_{t}\right) \right\vert +\left\vert
\varphi(Y_{t})-\varphi\left( u_{0}\right) \right\vert \right]
dt+\left\langle Y_{t}-u_{0},F\left( t,Y_{t},Z_{t}\right)
dt-dK_{t}\right\rangle \medskip \\
\leq\left[ R_{0}F_{u_{0},R_{0}}^{\#}\left( t\right) +R_{0}^{2}\gamma _{t}^{+}%
\right] dt+\left[ F_{u_{0},R_{0}}^{\#}\left( t\right) +2R_{0}\left\vert
\gamma_{t}\right\vert +2\left\vert \hat{u}_{0}\right\vert \right] \left\vert
Y_{t}-u_{0}\right\vert dt\medskip \\
\;\;+\left\vert Y_{t}-u_{0}\right\vert ^{2}dV_{t}+\dfrac{n_{p}}{2a}%
\left\vert Z_{t}\right\vert ^{2}%
\end{array}
\label{ineq 3 prop1}
\end{equation}
For $R_{0}=0$, inequality (\ref{ineq prop 1}) clearly follows from (\ref%
{ineq 3 prop1}) applying Proposition \ref{prop 1 Appendix} from Appendix.

For $R_{0}>0$ we moreover deduce, using once again Proposition \ref{prop 1
Appendix}, inequality (\ref{ineq 2 prop 1}).\hfill
\end{proof}

\begin{remark}
\label{remark 1}Denoting%
\begin{equation*}
\Theta=e^{V_{T}}\left\vert \eta-u_{0}\right\vert +\int
_{0}^{T}e^{V_{s}}\left\vert \hat{u}_{0}\right\vert ds+\int
_{0}^{T}e^{V_{s}}\left\vert F\left( s,u_{0},0\right) \right\vert ds
\end{equation*}
we deduce that, for all $t\in\left[ 0,T\right] :$%
\begin{equation}
\left\vert Y_{t}\right\vert \leq\left\vert u_{0}\right\vert
+C_{a,p}^{1/p}~e^{-V_{t}}\left( \mathbb{E}^{\mathcal{F}_{t}}\Theta^{p}%
\right) ^{1/p},\;\text{a.s.}  \label{ineq remark 1}
\end{equation}
\end{remark}

\begin{corollary}
\label{coroll 1}Let $p\geq 2$. We suppose moreover that there exist $%
r_{0},c_{0}>0$ such that%
\begin{equation*}
\varphi _{u_{0},r_{0}}^{\#}\overset{def}{=}\sup \left\{ \varphi \left(
u_{0}+r_{0}v\right) :\left\vert v\right\vert \leq 1\right\} \leq c_{0}~.
\end{equation*}%
Then%
\begin{equation}
\begin{array}{l}
\displaystyle r_{0}^{p/2}\mathbb{E}^{\mathcal{F}_{t}}\left(
\int_{t}^{T}e^{2V_{s}}d\left\updownarrow K\right\updownarrow _{s}\right)
^{p/2}\medskip \\
\displaystyle\leq C_{a,p}\mathbb{E}^{\mathcal{F}_{t}}\bigg[%
e^{pV_{T}}\left\vert \eta -u_{0}\right\vert ^{p}+\left( \varphi
_{u_{0},r_{0}}^{\#}-\varphi \left( u_{0}\right) \right) \left(
\int_{t}^{T}e^{2V_{s}}ds\right) ^{p/2}\medskip \\
\displaystyle\;\;+\left( \int_{t}^{T}e^{V_{s}}\left\vert \hat{u}%
_{0}\right\vert ds\right) ^{p}+\left( \int_{t}^{T}e^{V_{s}}\left\vert
F\left( s,u_{0},0\right) \right\vert ds\right) ^{p}\bigg]%
\end{array}
\label{ineq coroll 1}
\end{equation}
\end{corollary}

\begin{proof}
Let an arbitrary function $v\in C\left( \left[ 0,T\right] ;\mathbb{R}%
^{m}\right) $ such that $\left\Vert v\right\Vert _{T}\leq1$. From the
subdifferential inequality%
\begin{equation*}
\left\langle u_{0}+r_{0}v\left( t\right) -Y_{t},dK_{t}\right\rangle
+\varphi(Y_{t})dt\leq\varphi\left( u_{0}+r_{0}v\left( t\right) \right) dt,
\end{equation*}
we deduce that%
\begin{equation*}
r_{0}d\left\updownarrow K\right\updownarrow _{t}+\varphi(Y_{t})dt\leq
\left\langle Y_{t}-u_{0},dK_{t}\right\rangle +\varphi_{u_{0},r_{0}}^{\#}dt.
\end{equation*}
Since%
\begin{equation*}
\left\langle Y_{t}-u_{0},\hat{u}_{0}\right\rangle +\varphi\left(
u_{0}\right) \leq\varphi(Y_{t}),
\end{equation*}
then%
\begin{equation*}
r_{0}d\left\updownarrow K\right\updownarrow _{t}\leq\left\langle
Y_{t}-u_{0},dK_{t}\right\rangle +\left\vert \hat{u}_{0}\right\vert
\left\vert Y_{t}-u_{0}\right\vert dt+\left[ \varphi_{u_{0},r_{0}}^{\#}-%
\varphi\left( u_{0}\right) \right] dt
\end{equation*}
Therefore%
\begin{equation*}
\begin{array}{l}
r_{0}d\left\updownarrow K\right\updownarrow _{t}+\left\langle
Y_{t}-u_{0},F\left( t,Y_{t},Z_{t}\right) dt-dK_{t}\right\rangle \medskip \\
\leq\left( \varphi_{u_{0},r_{0}}^{\#}-\varphi\left( u_{0}\right) \right)
dt+\left\vert Y_{t}-u_{0}\right\vert \left( \left\vert \hat{u}%
_{0}\right\vert +\left\vert F\left( t,u_{0},0\right) \right\vert \right)
dt\medskip \\
+\left\vert Y_{t}-u_{0}\right\vert ^{2}dV_{t}+\dfrac{n_{p}}{2a}\left\vert
Z_{t}\right\vert ^{2}dt.%
\end{array}%
\end{equation*}
The inequality (\ref{ineq coroll 1}) follows using Proposition \ref{prop 1
Appendix}.\hfill
\end{proof}

\begin{proposition}[Uniqueness]
\label{prop 2 uniq}Let assumptions \textrm{(A}$_{1}-$\textrm{A}$_{3}$\textrm{%
)} be satisfied. Let $a,p>1$. If $\left( Y,Z\right) ,(\tilde{Y},\tilde{Z}%
)\in S_{m}^{0}\left[ 0,T\right] \times \Lambda _{m\times k}^{0}\left(
0,T\right) $ are two solutions of BSDE (\ref{BSVI}) corres\-ponding
respectively to $\eta $ and $\tilde{\eta}$ such that%
\begin{equation*}
\mathbb{E}\sup\limits_{s\in \left[ 0,T\right] }e^{pV_{s}}|Y_{s}-\tilde{Y}%
_{s}|^{p}<\infty ,
\end{equation*}%
then for all $t\in \left[ 0,T\right] ,$%
\begin{equation*}
e^{pV_{t}}|Y_{s}-\tilde{Y}_{s}|^{p}\leq \mathbb{E}^{\mathcal{F}_{t}}\left(
e^{pV_{T}}\left\vert \eta -\tilde{\eta}\right\vert ^{p}\right) ,\;\;\mathbb{P%
}\text{-a.s.}
\end{equation*}%
and there exists a constant $C_{a,p}$ such that $\mathbb{P}$-a.s.,\ for all $%
t\in \left[ 0,T\right] :$%
\begin{equation}
\begin{array}{r}
\displaystyle\mathbb{E}^{\mathcal{F}_{t}}\left[ \sup\limits_{s\in \left[ t,T%
\right] }e^{pV_{s}}|Y_{s}-\tilde{Y}_{s}|^{p}+\bigg(%
\int_{t}^{T}e^{2V_{s}}|Z_{s}-\tilde{Z}_{s}|^{2}ds\bigg)^{p/2}\right] \medskip
\\
\leq C_{a,p}\mathbb{E}^{\mathcal{F}_{t}}e^{pV_{T}}\left\vert \eta -\tilde{%
\eta}\right\vert ^{p}.%
\end{array}
\label{ineq uniq}
\end{equation}%
Moreover, the uniqueness of solution $\left( Y,Z\right) $ of BSDE (\ref{BSVI}%
) holds in $S_{m}^{1^{+},p}\left[ 0,T\right] \times \Lambda _{m\times
k}^{0}\left( 0,T\right) $.
\end{proposition}

\begin{proof}
Let $\left( Y,Z\right) $, $(\tilde{Y},\tilde{Z})\in S_{m}^{0}\left[ 0,T%
\right] \times\Lambda_{m\times k}^{0}\left( 0,T\right) $ be two solutions
corresponding to $\eta$ and $\tilde{\eta}$ respectively. Then there exists $%
p>1$ such that $Y,\tilde{Y}\in S_{m}^{p}\left[ 0,T\right] $ and
\begin{equation*}
Y_{t}-\tilde{Y}_{t}=\eta-\tilde{\eta}+\int_{t}^{T}dL_{s}-\int_{t}^{T}(Z_{s}-%
\tilde{Z}_{s})dB_{s}
\end{equation*}
where%
\begin{equation*}
L_{t}=\int _{0}^{t}\left[ (F\left( s,Y_{s},Z_{s}\right) -F\left( s,\tilde{Y}%
_{s},\tilde {Z}_{s}\right) )ds-(dK_{s}-d\tilde{K}_{s})\right] .
\end{equation*}
Since by (\ref{monoton}) $\langle Y_{s}-\tilde{Y}_{s},dK_{s}-d\tilde{K}%
_{s}\rangle\geq0,$ then, for all $a>1$,%
\begin{equation*}
\begin{array}{l}
\langle Y_{t}-\tilde{Y}_{t},dL_{t}\rangle\leq|Y_{t}-\tilde{Y}_{t}|^{2}\mu
_{t}dt+|Y_{t}-\tilde{Y}_{t}||Z_{t}-\tilde{Z}_{t}|\ell_{t}dt\medskip \\
\leq|Y_{t}-\tilde{Y}_{t}|^{2}\left( \mu_{t}+\dfrac{a}{2n_{p}}%
\ell_{t}^{2}\right) dt+\dfrac{n_{p}}{2a}|Z_{t}-\tilde{Z}_{t}|^{2}dt.%
\end{array}%
\end{equation*}
By Proposition \ref{prop 1 Appendix}, from Appendix, inequality (\ref{ineq
uniq}) follows.

Let now $p>1$ be such that $\left( Y,Z\right) ,(\tilde{Y},\tilde{Z})\in
S_{m}^{1^{+},p}\left[ 0,T\right] \times\Lambda_{m\times k}^{0}\left(
0,T\right) $ are two solutions of BSDE (\ref{BSVI}) corresponding
respectively to $\eta$ and $\tilde{\eta}$. From the definition of space $%
S_{m}^{1+,p}\left[ 0,T\right] $ there exists $a>1$ such that%
\begin{equation*}
\mathbb{E}\sup\limits_{t\in\left[ 0,T\right] }|e^{V_{t}^{a,p}}Y_{t}|^{p}<%
\infty,\;\;\mathbb{E}\sup\limits_{t\in\left[ 0,T\right] }|e^{V_{t}^{a,p}}%
\tilde{Y}_{t}|^{p}<\infty.
\end{equation*}
Consequently estimate (\ref{ineq uniq}) follows and uniqueness too.\hfill
\end{proof}

\section{BSVI - an existence result}

Using Proposition \ref{prop 1} we can prove now the existence of a triple $%
\left( Y,Z,K\right) $ which is a solution, in the sense of Definition \ref%
{def 1 - BSVI}, for BSVI (\ref{BSVI}). In order to obtain the absolute
continuity with respect to $dt$ for the process $K$ it is necessary to
impose a supplementary assumption.

Let $\left( u_{0},\hat{u}_{0}\right) \in\partial\varphi$ be fixed and%
\begin{equation}
\Theta_{T}^{a,p}\overset{def}{=}C_{a,p}e^{2p\left\Vert V\right\Vert _{T}}%
\left[ \left\vert \eta-u_{0}\right\vert ^{p}+\left( \int _{0}^{T}\left\vert
\hat{u}_{0}\right\vert ds\right) ^{p}+\left( \int _{0}^{T}\left\vert F\left(
s,u_{0},0\right) \right\vert ds\right) ^{p}\right] ,  \label{def of teta}
\end{equation}
where $a,p>1,$ $C_{a,p}$ is the constant given by Proposition \ref{prop 1}
and $V_{t}^{a,p}$ is defined by (\ref{def V}).

If there exists a constant $M$ such that%
\begin{equation*}
\left\vert \eta\right\vert +\int _{0}^{T}\left\vert F\left( s,u_{0},0\right)
\right\vert ds\leq M,\text{ a.s.}
\end{equation*}
then%
\begin{equation*}
\Theta_{T}^{a,p}\leq C_{a,p}e^{2p\left\Vert V\right\Vert _{T}}\left[ \left(
M+\left\vert u_{0}\right\vert \right) ^{p}+\left\vert \hat{u}_{0}\right\vert
^{p}T^{p}\right]
\end{equation*}
and by (\ref{ineq remark 1})%
\begin{equation*}
\left\vert Y_{t}\right\vert \leq\left\vert u_{0}\right\vert +\left( \mathbb{E%
}^{\mathcal{F}_{t}}\Theta_{T}^{a,p}\right) ^{1/p}\leq\left\vert
u_{0}\right\vert +C_{a,p}^{1/p}~e^{2\left\Vert V\right\Vert _{T}}\left[
M+\left\vert u_{0}\right\vert +\left\vert \hat{u}_{0}\right\vert T\right] ,%
\text{ a.s.}
\end{equation*}
We will make the following assumptions:

\begin{itemize}
\item[\textrm{(A}$_{4}$\textrm{)}] There exist $p\geq2,$ a positive
stochastic process $\beta\in L^{1}\left( \Omega\times\left( 0,T\right)
\right) $, a positive function $b\in L^{1}\left( 0,T\right) $ and a real
number $\kappa\geq0$, such that%
\begin{equation*}
\begin{array}{cl}
\left( i\right) & \mathbb{E}\varphi^{+}\left( \eta\right) <\infty ,\medskip
\\
\left( ii\right) & \text{for all\ }\left( u,\hat{u}\right) \in
\partial\varphi\mathit{\ }\text{and\ }z\in\mathbb{R}^{m\times k}:\medskip \\
& \;\;\left\langle \hat{u},F\left( t,u,z\right) \right\rangle \leq\dfrac {1}{%
2}\left\vert \hat{u}\right\vert ^{2}+\beta_{t}+b\left( t\right) \left\vert
u\right\vert ^{p}+\kappa\left\vert z\right\vert ^{2}\medskip \\
& \;\;\;\;\;\;\;\;\;\;\;\;\;\;\;\;d\mathbb{P}\otimes dt\text{-a.e., }\left(
\omega,t\right) \in\Omega\times\left[ 0,T\right] ,%
\end{array}%
\end{equation*}
and

\item[\textrm{(A}$_{5}$\textrm{)}] There exist $M,L>0$ and $\left( u_{0},%
\hat{u}_{0}\right) \in \partial \varphi $ such that:%
\begin{equation*}
\begin{array}{cl}
\left( i\right) & \mathbb{E}\varphi ^{+}\left( \eta \right) <\infty ,\medskip
\\
\left( ii\right) & \ell _{t}\leq L,\text{\ a.e., }t\in \left[ 0,T\right]
,\medskip \\
\left( iii\right) & \displaystyle\left\vert \eta \right\vert
+\int_{0}^{T}\left\vert F\left( s,u_{0},0\right) \right\vert ds\leq M,\text{%
\ a.s., }\omega \in \Omega ,\medskip \\
\left( iv\right) & \exists R_{0}\geq \left\vert u_{0}\right\vert
+C_{a,p}^{1/p}~e^{2\left\Vert V\right\Vert _{T}}\left[ M+\left\vert
u_{0}\right\vert +\left\vert \hat{u}_{0}\right\vert T\right] \medskip \\
& \displaystyle\text{\ \ \ \ \ such that }\mathbb{E}\int_{0}^{T}\left(
F_{R_{0}}^{\#}(s)\right) ^{2}ds<\infty%
\end{array}%
\end{equation*}
\end{itemize}

We note that, if $\left\langle \hat{u},F\left( t,u,z\right) \right\rangle
\leq0$, for all $\left( u,\hat{u}\right) \in\partial\varphi$, then condition
\textrm{(A}$_{4}$-$ii$\textrm{)} is satisfied with $\beta_{t}=b\left(
t\right) =\kappa=0$. For example, if $\varphi=I_{\bar{D}}$ (the convex
indicator of closed convex set $\bar{D}$) and $\mathbf{n}_{y}$ denotes the
unit outward normal vector to $\bar{D}$ at $y\in Bd\left( \bar{D}\right) $,
then condition $\left\langle n_{y},F\left( t,y,z\right) \right\rangle \leq0$
for all $y\in Bd\left( \bar{D}\right) $ yields \textrm{(A}$_{4}$-$ii$\textrm{%
)} with $\beta_{t}=b\left( t\right) =\kappa=0$. In this last case the It\^{o}%
's formula for $\psi\left( y\right) =\left[ dist_{\bar{D}}\left( y\right) %
\right] ^{2}$ and the uniqueness yields $K=0$.

We also remark that if $F\left( t,y,z\right) =F\left( y,z\right) $ then
assumptions \textrm{(A}$_{5}$\textrm{)} becomes%
\begin{equation*}
\left\vert \eta\right\vert +\mathbb{E}\varphi^{+}\left( \eta\right) \leq M,\;%
\text{a.s., }\omega\in\Omega.
\end{equation*}

\begin{theorem}[Existence]
\label{theor 1}Let $p\geq 2$ and assumptions \textrm{(A}$_{1}-$\textrm{A}$%
_{3}$\textrm{)} be satisfied with $s\rightarrow \mu _{s}=\mu \left( s\right)
$ and $s\rightarrow \ell _{s}=\ell \left( s\right) $ deterministic
processes. Suppose moreover that, for all $\rho \geq 0$,%
\begin{equation*}
\mathbb{E}\left\vert \eta \right\vert ^{p}+\mathbb{E}\left(
\int_{0}^{T}F_{\rho }^{\#}(s)ds\right) ^{p}<\infty ,
\end{equation*}%
and one of assumptions \textrm{(A}$_{4}$\textrm{) }or \textrm{(A}$_{5}$%
\textrm{)} is satisfied. Then there exists a unique pair $\left( Y,Z\right)
\in S_{m}^{p}\left[ 0,T\right] \times \Lambda _{m\times k}^{p}\left(
0,T\right) $ and a unique stochastic process $U\in \Lambda _{m}^{2}\left(
0,T\right) $ such that%
\begin{equation*}
\begin{array}{cl}
\left( a\right) & \displaystyle\int_{0}^{T}\left\vert
F(t,Y_{t},Z_{t})\right\vert dt<\infty ,\;\mathbb{P}\text{-a.s.},\medskip \\
\left( b\right) & Y_{t}\left( \omega \right) \in \mathrm{Dom}\left( \partial
\varphi \right) ,\;d\mathbb{P}\otimes dt\text{-\ a.e.\ }\left( \omega
,t\right) \in \Omega \times \left[ 0,T\right] ,\medskip \\
\left( c\right) & U_{t}\left( \omega \right) \in \partial \varphi \left(
Y_{t}\left( \omega \right) \right) ,\;d\mathbb{P}\otimes dt\text{\ -\ a.e.\ }%
\left( \omega ,t\right) \in \Omega \times \left[ 0,T\right]%
\end{array}%
\end{equation*}%
and for all $t\in \left[ 0,T\right] :$%
\begin{equation}
Y_{t}+\int_{t}^{T}U_{s}ds=\eta +\int_{t}^{T}F\left( s,Y_{s},Z_{s}\right)
ds-\int_{t}^{T}Z_{s}dB_{s},\text{\ a.s.}  \label{BSVI explicit form 2}
\end{equation}%
Moreover, uniqueness holds in $S_{m}^{1^{+}}\left[ 0,T\right] \times \Lambda
_{m\times k}^{0}\left( 0,T\right) ,$ where%
\begin{equation*}
S_{m}^{1^{+}}\left[ 0,T\right] \overset{def}{=}\bigcup_{p>1}S_{m}^{p}\left[
0,T\right] .
\end{equation*}
\end{theorem}

\begin{proof}
Let $\left( Y,Z\right) $, $(\tilde{Y},\tilde{Z})\in S_{m}^{1^{+}}\left[ 0,T%
\right] \times\Lambda_{m\times k}^{0}\left( 0,T\right) $ be two solutions.
Then $\exists p_{1},p_{2}>1$ such that $Y\in S_{m}^{p_{1}}\left[ 0,T\right] $%
, $\tilde{Y}\in S_{m}^{p_{2}}\left[ 0,T\right] $ and it follows that $Y,%
\tilde{Y}\in S_{m}^{p}\left[ 0,T\right] $, where $p=p_{1}\wedge p_{2}$.
Applying Proposition \ref{prop 2 uniq} we obtain the uniqueness.

To prove existence of a solution we can assume, without loss of generality,
that there exists $u_{0}\in \mathrm{Dom}\left( \varphi \right) $ such that%
\begin{equation}
0=\varphi \left( u_{0}\right) \leq \varphi \left( y\right) ,\;\forall y\in
\mathbb{R}^{m},  \label{pozitivity}
\end{equation}%
hence $0\in \partial \varphi \left( u_{0}\right) $, since, in the sense of
Definition \ref{def 1 - BSVI}, we can replace BSVI (\ref{BSVI}) by%
\begin{equation*}
\left\{
\begin{array}{l}
-dY_{t}+\partial \tilde{\varphi}\left( Y_{t}\right) dt\ni \tilde{F}\left(
t,Y_{t},Z_{t}\right) dt-Z_{t}dB_{t},\;0\leq t<T\medskip \\
Y_{T}=\eta ,%
\end{array}%
\,\right.
\end{equation*}%
where, for $\left( u_{0},\hat{u}_{0}\right) \in \partial \varphi $ fixed,%
\begin{equation*}
\begin{array}{l}
\tilde{\varphi}(y)\overset{def}{=}\varphi (y)-\varphi (u_{0})-\left\langle
\hat{u}_{0},y-u_{0}\right\rangle ,\;y\in \mathbb{R}^{d}\medskip \\
\tilde{F}(t,y,z)\overset{def}{=}F(t,y,z)-\hat{u}_{0},\;y\in \mathbb{R}%
^{d},\;t\in \left[ 0,T\right] .%
\end{array}%
\end{equation*}%
\textit{Step 1. Approximating problem.}

Let $\varepsilon\in(0,1]$ and the approximating equation%
\begin{equation}
Y_{t}^{\varepsilon}+\int _{t}^{T}\nabla\varphi_{\varepsilon}\left(
Y_{s}^{\varepsilon}\right) ds=\eta+\int _{t}^{T}F\left(
s,Y_{s}^{\varepsilon},Z_{s}^{\varepsilon}\right) ds-\int
_{t}^{T}Z_{s}^{\varepsilon}dB_{s},\text{\ a.s., }t\in\left[ 0,T\right] ,
\label{BSDE approx}
\end{equation}
$\nabla\varphi_{\varepsilon}$ is the gradient of the Yosida's regularization
$\varphi_{\varepsilon}$ of the function $\varphi$.

Using (\ref{pozitivity}) we obtain%
\begin{equation}
0=\varphi\left( u_{0}\right) \leq\varphi(J_{\varepsilon}y)\leq
\varphi_{\varepsilon}(y)\leq\varphi(y),\;J_{\varepsilon}\left( u_{0}\right)
=u_{0},\;\nabla\varphi_{\varepsilon}(u_{0})=0.  \label{ineq Yosida 2}
\end{equation}
It follows from \cite{br-de-hu-pa/03}, Theorem 4.2 (see also \cite{pa-ra/08}%
, Chapter 5) that equation (\ref{BSDE approx}) has an unique solution$%
\;\left( Y^{\varepsilon},Z^{\varepsilon}\right) \in S_{m}^{p}\left[ 0,T%
\right] \times\Lambda_{m\times k}^{p}\left( 0,T\right) $.$\medskip$\newline
\textit{Step 2. Boundedness of }$Y^{\varepsilon}$ \textit{and }$%
Z^{\varepsilon }$\textit{, without supplementary assumptions }\textrm{(A}$%
_{4}$\textrm{) }\textit{or} \textrm{(A}$_{5}$\textrm{)}\textit{.}

From Proposition \ref{prop 1}, applied for (\ref{BSDE approx}), we obtain,
for all $a>1$,%
\begin{equation}
\begin{array}{l}
\displaystyle\mathbb{E}^{\mathcal{F}_{t}}\left[ \sup\limits_{s\in \left[ t,T%
\right] }\left\vert e^{V_{s}}\left( Y_{s}^{\varepsilon }-u_{0}\right)
\right\vert ^{p}+\left( \int_{t}^{T}e^{2V_{s}}\varphi _{\varepsilon
}(Y_{s}^{\varepsilon })ds\right) ^{p/2}+\left(
\int_{t}^{T}e^{2V_{s}}\left\vert Z_{s}^{\varepsilon }\right\vert
^{2}ds\right) ^{p/2}\right] \medskip \\
\displaystyle\leq C_{a,p}\mathbb{E}^{\mathcal{F}_{t}}\left[
e^{pV_{T}}\left\vert \eta -u_{0}\right\vert ^{p}+\left(
\int_{t}^{T}e^{V_{s}}\left\vert F\left( s,u_{0},0\right) \right\vert
ds\right) ^{p}\right] .%
\end{array}
\label{ineq 1 th 1}
\end{equation}%
In particular there exists a constant independent of $\varepsilon $ such that%
\begin{equation}
\begin{array}{cl}
\left( a\right) & \displaystyle\mathbb{E}\left\Vert Y^{\varepsilon
}\right\Vert _{T}^{2}\leq \left( \mathbb{E}\left\Vert Y^{\varepsilon
}\right\Vert _{T}^{p}\right) ^{2/p}\leq C,\medskip \\
\left( b\right) & \displaystyle\mathbb{E}\int_{0}^{T}\left\vert
Z_{s}^{\varepsilon }\right\vert ^{2}ds\leq \bigg[\mathbb{E}\Big(%
\int_{0}^{T}\left\vert Z_{s}^{\varepsilon }\right\vert ^{2}ds\Big)^{p/2}%
\bigg]^{2/p}\leq C.%
\end{array}
\label{ineq 4 th 1}
\end{equation}%
Moreover, from (\ref{ineq remark 1}) we obtain%
\begin{equation}
\left\vert Y_{t}^{\varepsilon }\right\vert \leq \left\vert u_{0}\right\vert
+\left( \mathbb{E}^{\mathcal{F}_{t}}\Theta _{T}^{a,p}\right) ^{1/p}
\label{ineq 2 th 1}
\end{equation}%
where $\Theta _{T}^{a,p}$ is given by (\ref{def of teta}) with $\hat{u}%
_{0}=0 $ (since $\nabla \varphi _{\varepsilon }(u_{0})=0$).

Throughout the proof we shall fix $a=2$ (and then $V_{t}$ defined by (\ref%
{def V}), with $n_{p}=1\wedge\left( p-1\right) =1$, becomes $V_{t}=\int
_{0}^{t}\left[ \mu\left( s\right) +\ell^{2}\left( s\right) \right] ds$)$%
\medskip$\newline
\textit{Step 3. Boundedness of }$\nabla\varphi
_{\varepsilon}(Y_{s}^{\varepsilon})$.

Using the following stochastic subdifferential inequality (for proof see
Proposition 2.2, \cite{pa-ra/98})%
\begin{equation*}
\varphi _{\varepsilon }(Y_{t}^{\varepsilon })+\int_{t}^{T}\left\langle
\nabla \varphi _{\varepsilon }(Y_{s}^{\varepsilon }),\,dY_{s}^{\varepsilon
}\right\rangle \leq \varphi _{\varepsilon }(Y_{T}^{\varepsilon })=\varphi
_{\varepsilon }(\eta )\leq \varphi (\eta ),
\end{equation*}%
we deduce that, for all $t\in \left[ 0,T\right] ,$%
\begin{equation}
\begin{array}{r}
\displaystyle\varphi _{\varepsilon }(Y_{t}^{\varepsilon
})+\int_{t}^{T}\left\vert \nabla \varphi _{\varepsilon }(Y_{s}^{\varepsilon
})\right\vert ^{2}ds\leq \varphi (\eta )+\int_{t}^{T}\left\langle \nabla
\varphi _{\varepsilon }(Y_{s}^{\varepsilon }),F\left( s,Y_{s}^{\varepsilon
},Z_{s}^{\varepsilon }\right) \right\rangle ds\medskip \\
\displaystyle-\int_{t}^{T}\left\langle \nabla \varphi _{\varepsilon
}(Y_{s}^{\varepsilon }),\,Z_{s}^{\varepsilon }dB_{s}\right\rangle .%
\end{array}
\label{Ito for convex fct}
\end{equation}%
Since%
\begin{equation*}
\begin{array}{l}
\displaystyle\mathbb{E}\left( \int_{0}^{T}|\nabla \varphi _{\varepsilon
}\left( Y_{s}^{\varepsilon }\right) |^{2}|Z_{s}^{\varepsilon }|^{2}ds\right)
^{1/2}\leq \dfrac{1}{\varepsilon }\mathbb{E}\left[ \Big(\sup\limits_{s\in %
\left[ 0,T\right] }|Y_{s}^{\varepsilon }|\Big)\Big(\int_{0}^{T}|Z_{s}^{%
\varepsilon }|^{2}ds\Big)^{1/2}\right] \medskip \\
\displaystyle\leq \dfrac{1}{\varepsilon ^{2}}\mathbb{E}\Big(%
\sup\limits_{s\in \left[ 0,T\right] }|Y_{s}^{\varepsilon }|^{2}\Big)+\mathbb{%
E}\left( \int_{0}^{T}|Z_{s}^{\varepsilon }|^{2}ds\right) <\infty ,%
\end{array}%
\end{equation*}%
then%
\begin{equation*}
\mathbb{E}\int_{t}^{T}\left\langle \nabla \varphi _{\varepsilon
}(Y_{s}^{\varepsilon }),Z_{s}^{\varepsilon }dB_{s}\right\rangle =0.
\end{equation*}%
Under assumption \textrm{(A}$_{4}$\textrm{)}, since $\nabla \varphi
_{\varepsilon }(Y_{s}^{\varepsilon })\in \partial \varphi \left(
J_{\varepsilon }\left( Y_{s}^{\varepsilon }\right) \right) $, then%
\begin{equation*}
\begin{array}{l}
\left\langle \nabla \varphi _{\varepsilon }(Y_{s}^{\varepsilon }),F\left(
s,Y_{s}^{\varepsilon },Z_{s}^{\varepsilon }\right) \right\rangle \medskip \\
=\dfrac{1}{\varepsilon }\left\langle Y_{s}^{\varepsilon }-J_{\varepsilon
}\left( Y_{s}^{\varepsilon }\right) ,F\left( s,Y_{s}^{\varepsilon
},Z_{s}^{\varepsilon }\right) -F\left( s,J_{\varepsilon }\left(
Y_{s}^{\varepsilon }\right) ,Z_{s}^{\varepsilon }\right) \right\rangle
+\left\langle \nabla \varphi _{\varepsilon }(Y_{s}^{\varepsilon }),F\left(
s,J_{\varepsilon }\left( Y_{s}^{\varepsilon }\right) ,Z_{s}^{\varepsilon
}\right) \right\rangle \medskip \\
\leq \dfrac{1}{\varepsilon }\mu ^{+}\left( s\right) \left\vert
Y_{s}^{\varepsilon }-J_{\varepsilon }\left( Y_{s}^{\varepsilon }\right)
\right\vert ^{2}+\dfrac{1}{2}\left\vert \nabla \varphi _{\varepsilon
}(Y_{s}^{\varepsilon })\right\vert ^{2}+\beta _{s}+b\left( s\right)
\left\vert J_{\varepsilon }\left( Y_{s}^{\varepsilon }\right) \right\vert
^{p}+\kappa \left\vert Z_{s}^{\varepsilon }\right\vert ^{2}.%
\end{array}%
\end{equation*}%
From (\ref{def Yosida regul}) and inequality%
\begin{equation*}
\left\vert J_{\varepsilon }\left( Y_{s}^{\varepsilon }\right) \right\vert
\leq \left\vert J_{\varepsilon }\left( Y_{s}^{\varepsilon }\right)
-J_{\varepsilon }\left( u_{0}\right) \right\vert +\left\vert
u_{0}\right\vert \leq \left\vert Y_{s}^{\varepsilon }-u_{0}\right\vert
+\left\vert u_{0}\right\vert
\end{equation*}%
we have, for all $t\in \left[ 0,T\right] ,$%
\begin{equation*}
\begin{array}{l}
\displaystyle\mathbb{E}\varphi _{\varepsilon }(Y_{t}^{\varepsilon })+\dfrac{1%
}{2}\mathbb{E}\int_{t}^{T}\left\vert \nabla \varphi _{\varepsilon
}(Y_{s}^{\varepsilon })\right\vert ^{2}ds\medskip \\
\displaystyle\leq \mathbb{E}\varphi (\eta )+2\int_{t}^{T}\mu ^{+}\left(
s\right) \mathbb{E}\varphi _{\varepsilon }(Y_{s}^{\varepsilon })ds+\mathbb{E}%
\int_{t}^{T}\left( \beta _{s}+b\left( s\right) (\left\vert
Y_{s}^{\varepsilon }-u_{0}\right\vert +\left\vert u_{0}\right\vert
)^{p}+\kappa \left\vert Z_{s}^{\varepsilon }\right\vert ^{2}\right) ds%
\end{array}%
\end{equation*}%
that yields, via estimate (\ref{ineq 1 th 1}) and the backward Gronwall's
inequality, that there exists a constant $C>0$ independent of $\varepsilon
\in (0,1]$ such that%
\begin{equation}
\begin{array}{ll}
\left( a\right) & \displaystyle\mathbb{E}\varphi _{\varepsilon
}(Y_{t}^{\varepsilon })+\mathbb{E}\int_{0}^{T}\left\vert \nabla \varphi
_{\varepsilon }(Y_{s}^{\varepsilon })\right\vert ^{2}ds\leq C,\medskip \\
\left( b\right) & \mathbb{E}\left\vert Y_{t}^{\varepsilon }-J_{\varepsilon
}\left( Y_{t}^{\varepsilon }\right) \right\vert ^{2}\leq C\varepsilon .%
\end{array}
\label{ineq 3 th 1}
\end{equation}%
If we suppose \textrm{(A}$_{5}$\textrm{)} then, from (\ref{ineq 2 th 1}), we
infer that%
\begin{equation}
\begin{array}{l}
\left\vert Y_{t}^{\varepsilon }\right\vert \leq \left\vert u_{0}\right\vert
+\left( \mathbb{E}^{\mathcal{F}_{t}}\Theta _{T}^{2,p}\right) ^{1/p}\leq
\left\vert u_{0}\right\vert +C_{2,p}^{1/p}~e^{2\left\Vert V\right\Vert _{T}}%
\left[ M+\left\vert u_{0}\right\vert +\left\vert \hat{u}_{0}\right\vert T%
\right] \overset{def}{=}R_{0}%
\end{array}
\label{ineq 5 th 1}
\end{equation}%
Now%
\begin{equation*}
\begin{array}{l}
\left\langle \nabla \varphi _{\varepsilon }(Y_{s}^{\varepsilon }),F\left(
s,Y_{s}^{\varepsilon },Z_{s}^{\varepsilon }\right) \right\rangle \medskip \\
=\left\langle \nabla \varphi _{\varepsilon }(Y_{s}^{\varepsilon }),F\left(
s,Y_{s}^{\varepsilon },0\right) \right\rangle +\left\langle \nabla \varphi
_{\varepsilon }(Y_{s}^{\varepsilon }),F\left( s,Y_{s}^{\varepsilon
},Z_{s}^{\varepsilon }\right) -F\left( s,Y_{s}^{\varepsilon },0\right)
\right\rangle \medskip \\
\leq \dfrac{1}{2}\left\vert \nabla \varphi _{\varepsilon
}(Y_{s}^{\varepsilon })\right\vert ^{2}+|F_{R_{0}}^{\#}\left( s\right)
|^{2}+L^{2}\left\vert Z_{s}^{\varepsilon }\right\vert ^{2}%
\end{array}%
\end{equation*}%
Hence from (\ref{Ito for convex fct}) it follows that, for all $t\in \left[
0,T\right] ,$%
\begin{equation}
\mathbb{E}\varphi (J_{\varepsilon }\left( Y_{t}^{\varepsilon }\right) )+%
\dfrac{1}{2}\mathbb{E}\int_{t}^{T}\left\vert \nabla \varphi _{\varepsilon
}(Y_{s}^{\varepsilon })\right\vert ^{2}ds\leq \mathbb{E}\left( \varphi (\eta
)+\int_{t}^{T}|F_{R_{0}}^{\#}\left( s\right)
|^{2}ds+L^{2}\int_{t}^{T}\left\vert Z_{s}^{\varepsilon }\right\vert
^{2}ds\right)  \label{ineq 6 th 1}
\end{equation}%
and from (\ref{ineq 4 th 1}) we obtain boundedness inequalities (\ref{ineq 3
th 1}).

\textit{Step 4. Cauchy sequence and convergence.}

Let $\varepsilon,\delta\in(0,1].$

We can write%
\begin{equation*}
Y_{t}^{\varepsilon }-Y_{t}^{\delta }=\int_{t}^{T}dK_{s}^{\varepsilon ,\delta
}-\int_{t}^{T}Z_{s}^{\varepsilon }dB_{s}
\end{equation*}%
where%
\begin{equation*}
K_{t}^{\varepsilon ,\delta }=\int_{0}^{t}\left[ F\left( s,Y_{s}^{\varepsilon
},Z_{s}^{\varepsilon }\right) -F\left( s,Y_{s}^{\delta },Z_{s}^{\delta
}\right) -\nabla \varphi _{\varepsilon }\left( Y_{s}^{\varepsilon }\right)
+\nabla \varphi _{\delta }\left( Y_{s}^{\delta }\right) \right] ds.
\end{equation*}%
Then%
\begin{equation*}
\langle Y_{t}^{\varepsilon }-Y_{t}^{\delta },dK_{t}^{\varepsilon ,\delta
}\rangle \leq (\varepsilon +\delta )\langle \nabla \varphi _{\varepsilon
}(Y_{t}^{\varepsilon }),\nabla \varphi _{\delta }(Y_{t}^{\delta })\rangle
dt+|Y_{t}^{\varepsilon }-Y_{t}^{\delta }|^{2}dV_{t}+\dfrac{1}{4}%
|Z_{t}^{\varepsilon }-Z_{t}^{\delta }|^{2}dt,
\end{equation*}%
and by Proposition \ref{prop 1 Appendix}, with $p=2,$%
\begin{equation*}
\begin{array}{l}
\displaystyle\mathbb{E}\sup\limits_{s\in \left[ 0,T\right] }\left\vert
Y_{s}^{\varepsilon }-Y_{s}^{\delta }\right\vert ^{2}+\mathbb{E}%
\int_{0}^{T}\left\vert Z_{s}^{\varepsilon }-Z_{s}^{\delta }\right\vert
^{2}ds\leq C\mathbb{E}\int_{0}^{T}(\varepsilon +\delta )\left\langle \nabla
\varphi _{\varepsilon }(Y_{s}^{\varepsilon }),\nabla \varphi _{\delta
}(Y_{s}^{\delta })\right\rangle ds\medskip \\
\displaystyle\leq \dfrac{1}{2}C(\varepsilon +\delta )\left[ \mathbb{E}%
\int_{0}^{T}\left\vert \nabla \varphi _{\varepsilon }(Y_{s}^{\varepsilon
})\right\vert ^{2}ds+\mathbb{E}\int_{0}^{T}\left\vert \nabla \varphi
_{\delta }(Y_{s}^{\delta })\right\vert ^{2}ds\right] \leq C^{\prime
}(\varepsilon +\delta ).%
\end{array}%
\end{equation*}%
Hence there exist $\left( Y,Z,U\right) \in S_{m}^{2}\left[ 0,T\right] \times
\Lambda _{m\times k}^{2}\left( 0,T\right) \times \Lambda _{m}^{2}\left(
0,T\right) $ and a sequence $\varepsilon _{n}\searrow 0$ such that%
\begin{equation*}
\begin{array}{l}
Y^{\varepsilon _{n}}\rightarrow Y,\ \text{in }S_{m}^{2}\left[ 0,T\right]
\text{ and a.s. in }C\left( \left[ 0,T\right] ;\mathbb{R}^{m}\right)
,\medskip \\
Z^{\varepsilon _{n}}\rightarrow Z,\ \text{in }\Lambda _{m\times k}^{2}\left(
0,T\right) \text{ and a.s. in }L^{2}\left( 0,T;\mathbb{R}^{m\times k}\right)
,\medskip \\
\nabla \varphi _{\varepsilon }(Y^{\varepsilon })\rightharpoonup U,\ \text{%
weakly in }\Lambda _{m}^{2}\left( 0,T\right) ,\medskip \\
J_{\varepsilon _{n}}\left( Y^{\varepsilon _{n}}\right) \rightarrow Y,\ \text{%
in }\Lambda _{m}^{2}\left( 0,T\right) \text{ and a.s. in }L^{2}\left( 0,T;%
\mathbb{R}^{m}\right) .%
\end{array}%
\end{equation*}%
Passing to limit in (\ref{BSDE approx}) we conclude that%
\begin{equation*}
Y_{t}+\int_{t}^{T}U_{s}ds=\eta +\int_{t}^{T}F\left( s,Y_{s},Z_{s}\right)
ds-\int_{t}^{T}Z_{s}dB_{s},\;\text{a.s.}
\end{equation*}%
Since $\nabla \varphi _{\varepsilon }(Y_{s}^{\varepsilon })\in \partial
\varphi \left( J_{\varepsilon }\left( Y_{s}^{\varepsilon }\right) \right) $
then for all $A\in \mathcal{F}$, $0\leq s\leq t\leq T$ and $v\in S_{m}^{2}%
\left[ 0,T\right] ,$
\begin{equation*}
\mathbb{E}\int_{s}^{t}\mathbf{1}_{A}\left\langle \nabla \varphi
_{\varepsilon }(Y_{r}^{\varepsilon }),v_{r}-Y_{r}^{\varepsilon
}\right\rangle dr+\mathbb{E}\int_{s}^{t}\mathbf{1}_{A}\varphi
(J_{\varepsilon }\left( Y_{r}^{\varepsilon }\right) )dr\leq \mathbb{E}%
\int_{s}^{t}\mathbf{1}_{A}\varphi (v_{r})dr
\end{equation*}%
Passing to $\liminf $ for $\varepsilon =\varepsilon _{n}\searrow 0$ in the
above inequality we obtain that $U_{s}\in \partial \varphi \left(
Y_{s}\right) $. Hence $\left( Y,Z,U\right) \in S_{m}^{p}\left[ 0,T\right]
\times \Lambda _{m\times k}^{p}\left( 0,T\right) \times \Lambda
_{m}^{2}\left( 0,T\right) $ and $\left( Y,Z,K\right) ,$ with $%
K_{t}=\int_{0}^{t}U_{s}ds$, is the solution of BSVI (\ref{BSVI}).\newline
\textit{Step 5. Remarks in case }\textrm{(A}$_{5}$\textrm{)}.

Passing to $\liminf $ for $\varepsilon =\varepsilon _{n}\searrow 0$ in (\ref%
{ineq 5 th 1}) and (\ref{ineq 6 th 1}) it follows, using assumptions \textrm{%
(A}$_{5}$\textrm{)}, that the solution also satisfies%
\begin{equation*}
\begin{array}{ll}
\left( a\right) & \left\vert Y_{t}\right\vert \leq R_{0},\text{ a.s. for all
}t\in \left[ 0,T\right] ,\medskip \\
\left( b\right) & \displaystyle\mathbb{E}\varphi (Y_{t})+\dfrac{1}{2}\mathbb{%
E}\int_{t}^{T}\left\vert U_{s}\right\vert ^{2}ds\leq \mathbb{E}\left(
\varphi (\eta )+\int_{0}^{T}\left\vert F_{R_{0}}^{\#}\left( s\right)
\right\vert ^{2}ds+L^{2}\int_{0}^{T}\left\vert Z_{s}\right\vert
^{2}ds\right) .%
\end{array}%
\end{equation*}%
The proof is completed now.\hfill
\end{proof}

\begin{remark}
The existence Theorem \ref{theor 1} is well adapted to the Hilbert spaces
since we do not impose an assumption of type%
\begin{equation*}
\mathrm{Int}\left( \mathrm{Dom}\left( \varphi \right) \right) \neq \emptyset
,
\end{equation*}%
which is very restrictive for the infinite dimensional spaces. In the
context of the Hilbert spaces Theorem \ref{theor 1} holds in the same form
and one can give, as examples, partial differential backward stochastic
variational inequalities (see \cite{pa-ra/99}).
\end{remark}

\section{BSVI - a general existence result}

We replace now assumptions \textrm{(A}$_{5}$\textrm{)} with $\mathrm{Int}%
\left( \mathrm{Dom}\left( \varphi \right) \right) \neq \emptyset $.

\begin{theorem}[Existence]
\label{theor 2}Let $p\geq 2$ and assumptions \textrm{(A}$_{1}-$\textrm{A}$%
_{3}$\textrm{)} be satisfied with $s\rightarrow \mu _{s}=\mu \left( s\right)
$ and $s\rightarrow \ell _{s}=\ell \left( s\right) $ deterministic
processes. We suppose moreover that%
\begin{equation*}
\mathrm{Int}\left( \mathrm{Dom}\left( \varphi \right) \right) \neq \emptyset
\end{equation*}%
and for all $\rho \geq 0$%
\begin{equation*}
\mathbb{E}\left\vert \eta \right\vert ^{p}+\mathbb{E}\left(
\int_{0}^{T}F_{\rho }^{\#}(s)ds\right) ^{p}<\infty .
\end{equation*}%
Then there exists a unique triple $\left( Y,Z,K\right) \in S_{m}^{p}\left[
0,T\right] \times \Lambda _{m\times k}^{p}\left( 0,T\right) \times
S_{m}^{p}\left( 0,T\right) ,$ $\mathbb{E}\left\updownarrow
K\right\updownarrow _{T}^{p/2}<\infty ,$ such that for all $t\in \left[ 0,T%
\right] :$%
\begin{equation}
\left\{
\begin{array}{l}
\displaystyle Y_{t}+K_{T}-K_{t}=\eta +\int_{t}^{T}F\left(
s,Y_{s},Z_{s}\right) ds-\int_{t}^{T}Z_{s}dB_{s},\text{\ a.s.},\medskip \\
dK_{t}\in \partial \varphi \left( Y_{t}\right) dt,\text{ a.s.},\medskip \\
Y_{T}=\eta ,\text{ a.s.,}%
\end{array}%
\,\right.  \label{BSVI 2}
\end{equation}%
which means that BSVI (\ref{BSVI}) has a unique solution, and moreover%
\begin{equation*}
\mathbb{E}\left\Vert Y\right\Vert _{T}^{p}+\mathbb{E}\left\Vert K\right\Vert
_{T}^{p}+\mathbb{E}\left\updownarrow K\right\updownarrow _{T}^{p/2}+\mathbb{E%
}\int_{0}^{T}\left\vert Z_{t}\right\vert ^{2}dt<\infty .
\end{equation*}
\end{theorem}

\begin{proof}
The uniqueness was proved in Proposition \ref{prop 2 uniq}.$\medskip $%
\newline
\textit{Step 1. Existence under supplementary assumption}%
\begin{equation}
\begin{array}{l}
\exists M>0,\;u_{0}\in \mathrm{Int}\left( \mathrm{Dom}\left( \partial
\varphi \right) \right) \text{ \textit{such that}}\medskip \\
\displaystyle\mathbb{E}\left\vert \varphi \left( \eta \right) \right\vert
+\left\vert \eta \right\vert +\int_{0}^{T}\left\vert F\left(
s,u_{0},0\right) \right\vert ds\leq M,\ \text{a.s.\ }\omega \in \Omega%
\end{array}
\label{assumpt th 2}
\end{equation}%
Let $R_{0}$ defined by (\ref{ineq 5 th 1}) and denote%
\begin{equation*}
\zeta _{t}=\ell \left( t\right) +F_{R_{0}}^{\#}\left( t\right)
\end{equation*}%
By Theorem \ref{theor 1} there exists a unique $\left(
Y^{n},Z^{n},U^{n}\right) \in S_{m}^{p}\left[ 0,T\right] \times \Lambda
_{m\times k}^{p}\left( 0,T\right) \times \Lambda _{m}^{2}\left( 0,T\right) $
such that $U_{s}^{n}\in \partial \varphi \left( Y_{s}^{n}\right) $ and for
all $t\in \left[ 0,T\right] :$%
\begin{equation}
Y_{t}^{n}+\int_{t}^{T}U_{s}^{n}ds=\eta +\int_{t}^{T}F\left(
s,Y_{s}^{n},Z_{s}^{n}\right) \mathbf{1}_{\zeta _{t}\leq
n}ds-\int_{t}^{T}Z_{s}^{n}dB_{s},\text{\ a.s.}  \label{BSDE approx
2}
\end{equation}%
Moreover%
\begin{equation}
\sup\limits_{s\in \left[ 0,T\right] }\left\vert Y_{s}^{n}\right\vert \leq
R_{0},\text{\ a.s.}  \label{ineq 1 th 2}
\end{equation}%
and%
\begin{equation}
\mathbb{E}\left( \int_{0}^{T}\left\vert \varphi (Y_{s}^{n})\right\vert
ds\right) ^{p/2}+\mathbb{E}\left( \int_{0}^{T}\left\vert
Z_{s}^{n}\right\vert ^{2}ds\right) ^{p/2}\leq C.  \label{ineq 2 th 2}
\end{equation}%
Let $q=p/2$, $n_{q}=1\wedge \left( q-1\right) $, $a=2$ and $V_{t}^{2,q}$
given by (\ref{def V}).

Since%
\begin{equation*}
\begin{array}{l}
\left\langle Y_{t}^{n}-Y_{t}^{n+l},\left( F\left(
t,Y_{t}^{n},Z_{t}^{n}\right) \mathbf{1}_{\zeta _{t}\leq n}-U_{t}^{n}-F\left(
t,Y_{t}^{n+l},Z_{t}^{n+l}\right) \mathbf{1}_{\zeta _{t}\leq
n+l}+U_{t}^{n+l}\right) \right\rangle dt\medskip \\
\leq \left\langle Y_{t}^{n}-Y_{t}^{n+l},F\left( t,Y_{t}^{n},Z_{t}^{n}\right)
\right\rangle \left( \mathbf{1}_{\zeta _{t}\leq n}-\mathbf{1}_{\zeta
_{t}\leq n+l}\right) dt+\left\vert Y_{t}^{n}-Y_{t}^{n+l}\right\vert
^{2}dV_{t}^{2,q}\medskip \\
+\dfrac{n_{q}}{4}\left\vert Z_{t}^{n}-Z_{t}^{n+l}\right\vert ^{2}dt,%
\end{array}%
\end{equation*}%
then by Proposition \ref{prop 1 Appendix}, from Appendix, (with $a=2$) there
exists a constant depending only on $p$, such that%
\begin{equation*}
\begin{array}{l}
\displaystyle\mathbb{E}\sup\limits_{s\in \left[ 0,T\right] }\left\vert
Y_{s}^{n}-Y_{s}^{n+l}\right\vert ^{p/2}+\mathbb{E}\left(
\int_{0}^{T}\left\vert Z_{s}^{n}-Z_{s}^{n+l}\right\vert ^{2}ds\right)
^{p/4}\medskip \\
\displaystyle\leq C_{p}~e^{p\left\Vert V^{2,q}\right\Vert _{T}}\mathbb{E}%
\left( \int_{0}^{T}\mathbf{1}_{\zeta _{s}\geq n}\left\vert F\left(
s,Y_{s}^{n},Z_{s}^{n}\right) \right\vert ds\right) ^{p/2}.%
\end{array}%
\end{equation*}%
But%
\begin{equation*}
\begin{array}{l}
\displaystyle\mathbb{E}\left( \int_{0}^{T}\mathbf{1}_{\zeta _{s}\geq
n}\left\vert F\left( s,Y_{s}^{n},Z_{s}^{n}\right) \right\vert ds\right)
^{p/2}\leq \mathbb{E}\left( \int_{0}^{T}\mathbf{1}_{\zeta _{s}\geq n}\left(
F_{R_{0}}^{\#}\left( s\right) +\ell \left( s\right) \left\vert
Z_{s}^{n}\right\vert \right) ds\right) ^{p/2}\medskip \\
\displaystyle\leq C_{p}^{\prime }\mathbb{E}\left( \int_{0}^{T}\mathbf{1}%
_{\zeta _{s}\geq n}F_{R_{0}}^{\#}\left( s\right) ds\right) ^{p/2}\medskip \\
\displaystyle\;\;+C_{p}^{\prime }\left[ \mathbb{E}\left( \int_{0}^{T}\mathbf{%
1}_{\zeta _{s}\geq n}\ell ^{2}\left( s\right) ds\right) ^{p/2}\right]
^{1/2}\cdot \left[ \mathbb{E}\left( \int_{0}^{T}\left\vert Z^{n}\left(
s\right) \right\vert ^{2}ds\right) ^{p/2}\right] ^{1/2}\medskip \\
\displaystyle\leq C_{p}^{\prime }\left[ \mathbb{E}\left( \int_{0}^{T}\mathbf{%
1}_{\zeta _{s}\geq n}F_{R_{0}}^{\#}\left( s\right) ds\right) ^{p}\right]
^{1/2}\medskip \\
\displaystyle\;\;+C_{p}^{\prime }~C^{1/2}\left[ \mathbb{E}\left( \int_{0}^{T}%
\mathbf{1}_{\zeta _{s}\geq n}\ell ^{2}\left( s\right) ds\right) ^{p}\right]
^{1/2}\rightarrow 0,\;\text{as\ }n\rightarrow \infty .%
\end{array}%
\end{equation*}%
Hence there exists a pair $\left( Y,Z\right) \in S_{m}^{p/2}\left[ 0,T\right]
\times \Lambda _{m\times k}^{p/2}\left( 0,T\right) $ such that, as $%
n\rightarrow \infty $%
\begin{equation*}
\left( Y^{n},Z^{n}\right) \rightarrow \left( Y,Z\right) \;\text{in }%
S_{m}^{p/2}\left[ 0,T\right] \times \Lambda _{m\times k}^{p/2}\left(
0,T\right)
\end{equation*}%
In particular $Y_{0}^{n}\rightarrow Y_{0}$ in $\mathbb{R}^{m}$ and from
equation (\ref{BSDE approx 2}) it follows that%
\begin{equation*}
K_{\cdot }^{n}=\int_{0}^{\cdot }U_{s}^{n}ds\rightarrow K,\ \text{in }%
S_{m}^{0}\left[ 0,T\right] .
\end{equation*}%
Now by (\ref{ineq coroll 1}) for $V_{t}=V_{t}^{2,p}$ we obtain%
\begin{equation*}
\begin{array}{l}
\displaystyle\mathbb{E}\left( \int_{0}^{T}\left\vert U_{t}^{n}\right\vert
dt\right) ^{p/2}=\mathbb{E}\left\updownarrow K^{n}\right\updownarrow
_{T}^{p/2}\medskip \\
\displaystyle\leq Ce^{2p\left\Vert V\right\Vert _{T}}\left[ 1+T+\mathbb{E}%
\left\vert \eta \right\vert ^{p}+\mathbb{E}\left( \int_{0}^{T}\left\vert
F\left( t,u_{0},0\right) \right\vert dt\right) ^{p}\right]%
\end{array}%
\end{equation*}%
with $C=C\left( p,u_{0},\hat{u}_{0},r_{0},\varphi \right) .$

Therefore%
\begin{equation*}
\mathbb{E}\left\updownarrow K\right\updownarrow _{T}^{p/2}\leq
Ce^{2p\left\Vert V\right\Vert _{T}}\left[ 1+T+\mathbb{E}\left\vert
\eta\right\vert ^{p}+\mathbb{E}\left( \int _{0}^{T}\left\vert F\left(
s,u_{0},0\right) \right\vert ds\right) ^{p}\right] .
\end{equation*}
Passing to $\liminf$ as $n\rightarrow\infty$, eventually on a subsequence,
we deduce from (\ref{ineq 1 th 1}) and (\ref{ineq 2 th 1}) that%
\begin{equation*}
\sup\limits_{s\in\left[ 0,T\right] }\left\vert Y_{s}\right\vert \leq R_{0},%
\text{\ a.s.}
\end{equation*}
and%
\begin{equation*}
\mathbb{E}\left( \int _{0}^{T}\left\vert \varphi(Y_{s})\right\vert ds\right)
^{p/2}+\mathbb{E}\left( \int _{0}^{T}\left\vert Z_{s}\right\vert
^{2}ds\right) ^{p/2}\leq C.
\end{equation*}
To show that $\left( Y,Z,K\right) $ is solution of BSDE (\ref{BSVI 2}) it
remains to show that $dK_{t}\in\partial\varphi\left( Y_{t}\right) \left(
dt\right) $. Applying Corollary \ref{coroll 1 Appendix} we obtain $%
dK_{t}\in\partial\varphi\left( Y_{t}\right) \left( dt\right) $, since $%
dK_{t}^{n}=U_{t}^{n}dt\in\partial\varphi\left( Y_{t}^{n}\right) dt$.\newline
\textit{Step 2. Existence without supplementary assumption (\ref{assumpt th
2}).}

Let $\left( u_{0},\hat{u}_{0}\right) \in \partial \varphi $ such that $%
u_{0}\in \mathrm{Int}\left( \mathrm{Dom}\left( \varphi \right) \right) $ and
$\overline{B}\left( u_{0},r_{0}\right) \subset \mathrm{Dom}\left( \varphi
\right) .$ Recall that
\begin{equation*}
\varphi _{u_{0},r_{0}}^{\#}\overset{def}{=}\sup \left\{ \varphi \left(
u_{0}+r_{0}v\right) :\left\vert v\right\vert \leq 1\right\} <\infty .
\end{equation*}%
We introduce
\begin{equation*}
\eta ^{n}=\eta \mathbf{1}_{\left[ 0,n\right] }\left( \left\vert \eta
\right\vert +\left\vert \varphi \left( \eta \right) \right\vert \right)
+u_{0}\mathbf{1}_{(n,\infty )}\left( \left\vert \eta \right\vert +\left\vert
\varphi \left( \eta \right) \right\vert \right)
\end{equation*}%
and%
\begin{equation*}
F^{n}\left( t,y,z\right) =F\left( s,y,z\right) -F\left( s,u_{0},0\right)
\mathbf{1}_{\left\vert F\left( s,u_{0},0\right) \right\vert \geq n}
\end{equation*}%
Clearly%
\begin{equation*}
\left\vert \eta ^{n}\right\vert +\left\vert \varphi \left( \eta _{n}\right)
\right\vert +\left\vert F^{n}\left( t,u_{0},0\right) \right\vert \leq
3n+\left\vert \varphi \left( u_{0}\right) \right\vert .
\end{equation*}%
By \textit{Step 1}, for each $n\in \mathbb{N}^{\ast }$ there exists a unique
triple $\left( Y^{n},Z^{n},K^{n}\right) \in S_{m}^{p}\left[ 0,T\right]
\times \Lambda _{m\times k}^{p}\left( 0,T\right) \times S_{m}^{p/2}\left(
0,T\right) $ solution of BSDE%
\begin{equation}
Y_{t}^{n}+\left( K_{T}^{n}-K_{t}^{n}\right) =\eta
^{n}+\int_{t}^{T}F^{n}\left( s,Y_{s}^{n},Z_{s}^{n}\right)
ds-\int_{t}^{T}Z_{s}^{n}dB_{s},\text{\ a.s.}  \label{BSDE approx 3}
\end{equation}%
From Corollary \ref{coroll 1} and Proposition \ref{prop 2 uniq} we infer
that there exists a constant $C_{p}$ such that%
\begin{equation}
\begin{array}{l}
\displaystyle\mathbb{E}r_{0}^{p/2}\left\updownarrow K^{n}\right\updownarrow
_{T}^{p/2}+\mathbb{E}\sup\limits_{s\in \left[ 0,T\right] }\left\vert \left(
Y_{s}^{n}-u_{0}\right) \right\vert ^{p}+\mathbb{E}\Big(\int_{0}^{T}\left%
\vert \varphi (Y_{s}^{n})-\varphi \left( u_{0}\right) \right\vert ds\Big)%
^{p/2}\medskip \\
\displaystyle\;\;\;+\mathbb{E}\Big(\int_{0}^{T}\left\vert
Z_{s}^{n}\right\vert ^{2}ds\Big)^{p/2}\leq C_{p}e^{2p\left\Vert V\right\Vert
_{T}}\bigg[\left[ \varphi _{u_{0},r_{0}}^{\#}-\varphi \left( u_{0}\right) %
\right] ^{p/2}T^{p/2}+\left\vert \hat{u}_{0}\right\vert ^{p}T^{p}\medskip \\
\displaystyle\;\;\;+\mathbb{E}\left\vert \eta ^{n}-u_{0}\right\vert ^{p}+%
\mathbb{E}\Big(\int_{0}^{T}\left\vert F^{n}\left( s,u_{0},0\right)
\right\vert ds\Big)^{p}\bigg]\medskip \\
\displaystyle\leq C_{p}e^{2p\left\Vert V\right\Vert _{T}}\bigg[\left[
\varphi _{u_{0},r_{0}}^{\#}-\varphi \left( u_{0}\right) \right]
^{p/2}T^{p/2}+\left\vert \hat{u}_{0}\right\vert ^{p}T^{p}+\mathbb{E}%
\left\vert \eta -u_{0}\right\vert ^{p}+\medskip \\
\displaystyle\;\;\;+\mathbb{E}\Big(\int_{0}^{T}\left\vert F\left(
s,u_{0},0\right) \right\vert ds\Big)^{p}\bigg]%
\end{array}
\label{ineq 3 th 2}
\end{equation}%
Remark that $p\geq 2$ is required only to obtain the estimate of $\mathbb{E}%
\left\updownarrow K^{n}\right\updownarrow _{T}^{p/2}$.

Since%
\begin{equation*}
\begin{array}{l}
\langle Y_{s}^{n}-Y_{s}^{n+l},F^{n}\left( s,Y_{s}^{n},Z_{s}^{n}\right)
-F^{n+l}\left( s,Y_{s}^{n+l},Z_{s}^{n+l}\right) \rangle \medskip \\
\displaystyle\leq \left\vert Y_{s}^{n}-Y_{s}^{n+l}\right\vert \left\vert
F\left( s,u_{0},0\right) \right\vert \mathbf{1}_{\left\vert F\left(
s,u_{0},0\right) \right\vert \geq n}~+\left\vert
Y_{s}^{n}-Y_{s}^{n+l}\right\vert ^{2}dV_{t}+\frac{1}{4}\left\vert
Z_{s}^{n}-Z_{s}^{n+l}\right\vert ^{2}ds%
\end{array}%
\end{equation*}%
then by Proposition \ref{prop 1 Appendix} we obtain%
\begin{equation*}
\begin{array}{l}
\displaystyle\mathbb{E}\Big(\sup\limits_{s\in \left[ 0,T\right] }\left\vert
Y_{s}^{n}-Y_{s}^{n+l}\right\vert ^{p}\Big)+\mathbb{E}\Big(%
\int_{0}^{T}\left\vert Z_{s}^{n}-Z_{s}^{n+l}\right\vert ^{2}ds\Big)%
^{p/2}\medskip \\
\displaystyle\leq C_{p}e^{2p\left\Vert V\right\Vert _{T}}\left[ \mathbb{E}%
\left( \left\vert \eta -u_{0}\right\vert ^{p}\mathbf{1}_{\left\vert \eta
\right\vert +\left\vert \varphi \left( \eta \right) \right\vert \geq
n}\right) +\mathbb{E}\Big(\int_{0}^{T}\left\vert F\left( s,u_{0},0\right)
\right\vert \mathbf{1}_{\left\vert F\left( s,u_{0},0\right) \right\vert \geq
n}\Big)^{p}\right] .%
\end{array}%
\end{equation*}%
Hence there exists a pair $\left( Y,Z\right) \in S_{m}^{p}\left[ 0,T\right]
\times \Lambda _{m\times k}^{p}\left( 0,T\right) $ such that
\begin{equation*}
\left( Y^{n},Z^{n}\right) \rightarrow \left( Y,Z\right) ,\;\text{as }%
n\rightarrow \infty ,\;\text{in }S_{m}^{p}\left[ 0,T\right] \times \Lambda
_{m\times k}^{p}\left( 0,T\right) .
\end{equation*}%
In particular $Y_{0}^{n}\rightarrow Y_{0}$ in $\mathbb{R}^{m}$. From
equation (\ref{BSDE approx 3}) we have%
\begin{equation*}
K^{n}\rightarrow K\;\text{in }S_{m}^{0}\left[ 0,T\right] ,
\end{equation*}%
and for all $t\in \left[ 0,T\right] $%
\begin{equation*}
Y_{t}+K_{T}-K_{t}=\eta +\int_{t}^{T}F\left( s,Y_{s},Z_{s}\right)
ds-\int_{t}^{T}Z_{s}dB_{s},\text{\ a.s.}
\end{equation*}%
Letting $n\rightarrow \infty $ and applying Proposition \ref{prop 2 Appendix}
we can assert that estimate (\ref{ineq 3 th 2}) holds without $n$. To
complete the proof remark that from $dK_{t}^{n}\in \partial \varphi \left(
Y_{t}^{n}\right) dt$ we can infer, using Corollary \ref{coroll 1 Appendix},
that $dK_{t}\in \partial \varphi \left( Y_{t}\right) dt$.

Therefore $\left( Y,Z,K\right) $ is solution of BSDE (\ref{BSVI 2}) in the
sense of Definition \ref{def 1 - BSVI}.\hfill
\end{proof}

\begin{remark}
When $\mu$ and $\ell$ are stochastic processes we obtain, with similar
proofs as in Theorems \ref{theor 1} and \ref{theor 2}, the existence of a
solution in the space%
\begin{equation*}
\mathbb{U}_{m,k}^{p}\left( 0,T\right) \overset{def}{=}\left\{ \left(
Y,Z\right) \in S_{m}^{0}\left[ 0,T\right] \times\Lambda_{m\times
k}^{0}\left( 0,T\right) :\left\Vert \left( Y,Z\right) \right\Vert
_{a,p}<\infty,\;\forall a>1\right\} ,
\end{equation*}
where%
\begin{equation*}
\left\Vert \left( Y,Z\right) \right\Vert _{a,p}^{p}\overset{def}{=}\mathbb{E}%
\Big(\sup\limits_{s\in\left[ 0,T\right] }e^{pV_{s}^{a,p}}|Y_{s}|^{p}\Big)+%
\mathbb{E}\Big(\int _{0}^{T}e^{2V_{s}^{a,p}}|Z_{s}|^{2}ds\Big)^{p/2}~.
\end{equation*}
\end{remark}

\section{Appendix}

In this section we first present some useful and general estimates on $%
\left( Y,Z\right) \in S_{m}^{0}\left[ 0,T\right] \times\Lambda_{m\times
k}^{0}\left( 0,T\right) $ satisfying an identity of type%
\begin{equation*}
Y_{t}=Y_{T}+\int_{t}^{T}dK_{s}-\int_{t}^{T}Z_{s}dB_{s},\;t\in\left[ 0,T%
\right] ,\;\mathbb{P}\text{-a.s.,}
\end{equation*}
where $K\in S_{m}^{0}\left[ 0,T\right] $ and $K_{\cdot}\left( \omega\right)
\in BV\left( \left[ 0,T\right] ;\mathbb{R}^{m}\right) $ $\mathbb{P}$-a.s., $%
\omega\in\Omega.$

The following results and their proofs are given in the monograph of E.
Pardoux, A. R\u{a}\c{s}canu \cite{pa-ra/08}, Annex C (a forthcoming
publication).

Assume there exist

$\diamondsuit\;D,R,N$ progressively measurable increasing continuous
stochastic processes with $D_{0}=R_{0}=N_{0}=0$,

$\diamondsuit\;V$ progressively measurable bounded-variation continuous
stochastic process with $V_{0}=0$,

$\diamondsuit\;a,p>1$,

such that, as signed measures on $\left[ 0,T\right] :$%
\begin{equation}
dD_{t}+\left\langle Y_{t},dK_{t}\right\rangle \leq\left( \mathbf{1}_{p\geq
2}dR_{t}+|Y_{t}|dN_{t}+|Y_{t}|^{2}dV_{t}\right) +\dfrac{n_{p}}{2a}\left\vert
Z_{t}\right\vert ^{2}dt  \label{ineq D,R,N}
\end{equation}
where $n_{p}=\left( p-1\right) \wedge1$.

Let $\left\Vert Ye^{V}\right\Vert _{\left[ t,T\right] }\overset{def}{=}%
\sup\limits_{s\in\left[ t,T\right] }\left\vert Y_{s}e^{V_{s}}\right\vert $
and $\left\Vert Ye^{V}\right\Vert _{T}\overset{def}{=}\left\Vert
Ye^{V}\right\Vert _{\left[ 0,T\right] }$ .

\begin{proposition}
\label{prop 1 Appendix}Assume (\ref{ineq D,R,N}) and that%
\begin{equation*}
\mathbb{E}\left\Vert Ye^{V}\right\Vert _{T}^{p}+\mathbb{E}\left(
\int_{0}^{T}e^{2V_{s}}\mathbf{1}_{p\geq 2}dR_{s}\right) ^{p/2}+\mathbb{E}%
\left( \int_{0}^{T}e^{V_{s}}dN_{s}\right) ^{p}<\infty .
\end{equation*}%
Then there exists a positive constant $C_{a,p}$, depending only of $a,p$,
such that, $\mathbb{P}$-a.s., for all $t\in \left[ 0,T\right] :$%
\begin{equation}
\begin{array}{l}
\displaystyle\mathbb{E}^{\mathcal{F}_{t}}\left[ \sup\limits_{s\in \left[ t,T%
\right] }\left\vert e^{V_{s}}Y_{s}\right\vert ^{p}+\left(
\int_{t}^{T}e^{2V_{s}}dD_{s}\right) ^{p/2}+\left(
\int_{t}^{T}e^{2V_{s}}\left\vert Z_{s}\right\vert ^{2}ds\right) ^{p/2}\right]
\medskip \\
\displaystyle\;\;+\mathbb{E}^{\mathcal{F}_{t}}\left[ \int_{t}^{T}e^{pV_{s}}%
\left\vert Y_{s}\right\vert ^{p-2}\mathbf{1}_{Y_{s}\neq
0}dD_{s}+\int_{t}^{T}e^{pV_{s}}\left\vert Y_{s}\right\vert ^{p-2}\mathbf{1}%
_{Y_{s}\neq 0}\left\vert Z_{s}\right\vert ^{2}ds\right] \medskip \\
\displaystyle\leq C_{a,p}\;\mathbb{E}^{\mathcal{F}_{t}}\left[ \left\vert
e^{V_{T}}Y_{T}\right\vert ^{p}+\left( \int_{t}^{T}e^{2V_{s}}\mathbf{1}%
_{p\geq 2}dR_{s}\right) ^{p/2}+\left( \int_{t}^{T}e^{V_{s}}dN_{s}\right) ^{p}%
\right] .%
\end{array}
\label{ineq 1 appendix}
\end{equation}%
In particular for all $t\in \left[ 0,T\right] :$%
\begin{equation*}
\left\vert Y_{t}\right\vert ^{p}\leq C_{a,p}\;\mathbb{E}^{\mathcal{F}_{t}}%
\left[ \left( \left\vert Y_{T}\right\vert ^{p}+\mathbf{1}_{p\geq
2}R_{T}^{p}+N_{T}^{p}\right) e^{p\left\Vert \left( V_{\cdot }-V_{t}\right)
^{+}\right\Vert _{\left[ t,T\right] }}\right] ,\;\mathbb{P}\text{-a.s.,}
\end{equation*}%
Moreover if there exists a constant $b\geq 0$ such that for all $t\in \left[
0,T\right] :$%
\begin{equation*}
\left\vert e^{V_{T}-V_{t}}Y_{T}\right\vert +\left( \int_{t}^{T}e^{2\left(
V_{s}-V_{t}\right) }\mathbf{1}_{p\geq 2}dR_{s}\right)
^{1/2}+\int_{t}^{T}e^{\left( V_{s}-V_{t}\right) }dN_{s}\leq b,\;\text{a.s.}
\end{equation*}%
then for all $t\in \left[ 0,T\right] :$%
\begin{equation}
\left\vert Y_{t}\right\vert ^{p}+\mathbb{E}^{\mathcal{F}_{t}}\left(
\int_{t}^{T}e^{2\left( V_{s}-V_{t}\right) }\left\vert Z_{s}\right\vert
^{2}ds\right) ^{p/2}\leq b^{p}C_{a,p},\;\;\mathbb{P}\text{-a.s.}
\label{ineq 2 appendix}
\end{equation}
\end{proposition}

The following results provides a criterion for passing to the limit in
Stieltjes integral (for the proofs we refer the reader to \cite{pa-ra/08},
Chapter I).

\begin{proposition}
\label{prop 2 Appendix}Let $Y,K,Y^{n},K^{n}$ be $C\left( \left[ 0,T\right] ;%
\mathbb{R}^{m}\right) $-valued random variables, $n\in\mathbb{N}$. Assume%
\begin{equation*}
\begin{array}{cl}
\left( i\right) & \exists p>0\;\text{such that }\sup\limits_{n\in \mathbb{N}%
^{\ast}}\mathbb{E}\left\updownarrow K^{n}\right\updownarrow _{T}^{p}<\infty%
\text{,}\medskip \\
\left( ii\right) & \left( \left\Vert Y^{n}-Y\right\Vert _{T}+\left\Vert
K^{n}-K\right\Vert _{T}\right)
\xrightarrow[]{prob.}%
0,\;\text{as\ }n\rightarrow\infty,\medskip \\
& \text{ \ \ i.e. }\forall\varepsilon>0,\;\mathbb{P}\left\{ \left(
\left\Vert Y^{n}-Y\right\Vert _{T}+\left\Vert K^{n}-K\right\Vert _{T}\right)
>\varepsilon\right\} \rightarrow0,\;\text{as\ }n\rightarrow\infty.%
\end{array}%
\end{equation*}
Then, for all $0\leq s\leq t\leq T:$%
\begin{equation*}
\int_{s}^{t}\left\langle Y_{r}^{n},dK_{r}^{n}\right\rangle
\xrightarrow[]{prob.}%
\int_{s}^{t}\left\langle Y_{r},dK_{r}\right\rangle ,\;\text{as}\
n\rightarrow \infty,
\end{equation*}
and moreover,%
\begin{equation*}
\mathbb{E}\left\updownarrow K\right\updownarrow _{T}^{p}\leq\liminf
_{n\rightarrow+\infty}\mathbb{E}\left\updownarrow K^{n}\right\updownarrow
_{T}^{p}\;.
\end{equation*}
\end{proposition}

\begin{corollary}
\label{coroll 1 Appendix}Let\ the assumptions of Proposition \ref{prop 2
Appendix} be satisfied. If $A:\mathbb{R}^{m}\rightrightarrows \mathbb{R}^{m}$
is a (multivalued) maximal monotone operator then the following implication
holds%
\begin{equation*}
dK_{t}^{n}\in A\left( Y_{t}^{n}\right) dt\;\text{on }\left[ 0,T\right] ,%
\text{\ a.s.}\;\Rightarrow\;dK_{t}\in A\left( Y_{t}\right) dt\;\text{on }%
\left[ 0,T\right] ,\text{\ a.s.}
\end{equation*}
In particular if $\varphi:\mathbb{R}^{d}\rightarrow]-\infty,+\infty]$ is a
proper convex l.s.c. function then%
\begin{equation*}
dK_{t}^{n}\in\partial\varphi\left( Y_{t}^{n}\right) dt\;\text{on }\left[ 0,T%
\right] ,\text{\ a.s.}\;\Rightarrow\;dK_{t}\in\partial\varphi\left(
Y_{t}\right) dt\;\text{on }\left[ 0,T\right] ,\text{\ a.s.}
\end{equation*}
\end{corollary}


\begin{thebibliography}{99}
\bibitem{br/73} H. Br\'{e}zis, \textit{Op\'{e}rateurs maximaux monotones et
semigroupes de contractions dans les espaces de Hilbert}, North-Holland,
Amsterdam, 1973.\vspace{-0.1cm}

\bibitem{br-de-hu-pa/03} Ph. Briand, B. Delyon, Y. Hu, E. Pardoux, L.
Stoica, $L^{p}$\textit{\ solutions of backward stochastic differential
equations}, Stochastic Processes and their Applications 108 (2003), 109-129.%
\vspace{-0.1cm}

\bibitem{cv-ka/96} J. Cvitanic, I. Karatzas, \textit{Backward stochastic
differential equations with reflection and Dynkin games}, Ann. Probab., 24
(1996), 2024-2056.\vspace{-0.1cm}

\bibitem{ka-ka-pa/97} N. El Karoui, C. Kapoudjian, E. Pardoux, S. Peng, M.C.
Quenez, \textit{Reflected solutions of backward SDE's and related obstacle
problems for PDE's}, Ann.Probab. 25 (1997) 702-737.\vspace{-0.1cm}

\bibitem{ge-pa/96} A. Gegout-Petit, E. Pardoux, \textit{Equations diff\'{e}%
rentielles stochastiques r\'{e}trogrades r\'{e}fl\'{e}chies dans un convexe}%
, Stochastics Stochastics Rep., 57 (1996), 111--128.\vspace{-0.1cm}

\bibitem{ma-ro/12} L. Maticiuc, E. Rotenstein, \textit{Numerical schemes for
multivalued backward stochastic differential systems}, Central European
Journal of Mathematics, 10 (2) (2012), 693--702.\vspace{-0.1cm}

\bibitem{pa/99} E. Pardoux, \textit{BSDEs, weak convergence and
homogenization of semilinear PDEs}, Nonlinear Analysis, Differential
Equations and Control (Montreal, QC, 1998), Kluwer Academic Publishers,
Dordrecht (1999), 503--549.\vspace{-0.1cm}

\bibitem{pa-pe/90} E. Pardoux, S. Peng, \textit{Adapted solution of a
backward stochastic differential equation}, System and control letters, 14
(1990), 55-61.\vspace{-0.1cm}

\bibitem{pa-pe/92} E. Pardoux, S. Peng, \textit{Backward SDE's and
quasilinear parabolic PDE's}, Stochastic PDE and Their Applications, LNCIS
176, Springer (1992) 200-217.\vspace{-0.1cm}

\bibitem{pa-ra/98} E. Pardoux, A. R\u{a}\c{s}canu, \textit{Backward
stochastic differential equations with subdifferential operator and related
variational inequalities}, Stochastic Processes and their Applications 76
(1998) 191-215.\vspace{-0.1cm}

\bibitem{pa-ra/99} E. Pardoux, A. R\u{a}\c{s}canu, \textit{Backward
stochastic variational inequalities}, Stochastics Stochastics Rep. 67 (3-4)
(1999) 159--167.\vspace{-0.1cm}

\bibitem{pa-ra/08} E. Pardoux, A. R\u{a}\c{s}canu, \textit{Stochastic
differential equations, Backward SDEs, Partial differential equations},
Stochastic Modelling and Applied Probability series, Springer, in press,
2014.
\end{thebibliography}
\end{document}